\documentclass{amsproc}
\usepackage{euscript}
\usepackage{cases}
\usepackage{mathrsfs}
\usepackage{bbm}
\usepackage{amssymb}
\usepackage{amsfonts,amsmath,amsxtra,mathdots,mathabx,pifont}
\usepackage{color}
\usepackage{hyperref}
\usepackage{tikz}
\usepackage{appendix,upgreek}

\allowdisplaybreaks

\DeclareFontFamily{U}{matha}{\hyphenchar\font45}
\DeclareFontShape{U}{matha}{m}{n}{
	<5> <6> <7> <8> <9> <10> gen * matha
	<10.95> matha10 <12> <14.4> <17.28> <20.74> <24.88> matha12
}{}
\DeclareSymbolFont{matha}{U}{matha}{m}{n}

\DeclareMathSymbol{\Lt}{3}{matha}{"CE}
\DeclareMathSymbol{\Gt}{3}{matha}{"CF}

\def\valpha{\text{\scalebox{0.88}[1.02]{$\alpha$}}}   
\def\vepsilon{\upvarepsilon} 
\def\vchi{\text{\raisebox{0.6 \depth}{\scalebox{0.9}[1.1]{$\chi$}}}} 
\def\vlambda{\text{\scalebox{0.9}[1]{$\lambda$}}}

\def\vnu{\text{{\scalebox{0.9}[1]{$\nu$}}}} 
\def\uppii{\text{\scalebox{0.8}[0.96]{$\uppi$}}}

\def\RC{\mathrm {C}}
\def\RN{\mathrm{N}}

\def\sumx{\sideset{}{^{\scriptscriptstyle\text{\ding{72}}}}\sum}
\def\sumst{\sideset{}{^{\star}}\sum}
\def\sumf{\sideset{}{^\flat}\sum} 
\def\sumsp{\sideset{}{^{ {\scriptscriptstyle  \text{\ding{73}} }    }}\sum}

\def\SD{\text{\raisebox{- 2 \depth}{\scalebox{1.1}{$ \text{\usefont{U}{BOONDOX-calo}{m}{n}D}   $}}}}

\def\SC{\text{\raisebox{- 2 \depth}{\scalebox{1.1}{$ \text{\usefont{U}{BOONDOX-calo}{m}{n}C}  $}}}}

\def\SE{\text{\raisebox{- 2 \depth}{\scalebox{1.1}{$ \text{\usefont{U}{BOONDOX-calo}{m}{n}E}  $}}}}

\def\SF{\text{\raisebox{- 2 \depth}{\scalebox{1.1}{$ \text{\usefont{U}{BOONDOX-calo}{m}{n}F}  $}}}}

\def\CaloO {\text{\raisebox{- 2 \depth}{\scalebox{1.1}{$ \text{\usefont{U}{BOONDOX-calo}{m}{n}O}  $}}}}

\def\SDH{\text{\raisebox{- 6 \depth}{\scalebox{1.06}{$ \text{\usefont{U}{dutchcal}{m}{n}H}  $}}}}

\def\shskip{\hspace{0.5pt}}

\def\frc{\mathfrak{c}}
\def\frm{\mathfrak{m}}
\def\frn{\mathfrak{n}}
\def\frp{\mathfrak{p}}

\newcommand{\delete}[1]{}

\theoremstyle{plain}

\newtheorem{thm}{Theorem} \newtheorem{cor}[thm]{Corollary}
\newtheorem{lem}{Lemma}[section]  \newtheorem{prop}[thm]{Proposition}

\theoremstyle{remark} 
\newtheorem{remark}{Remark}[section]

\numberwithin{equation}{section}

\begin{document}

	\title[Symmetric Square Large Sieve and Prime Geodesic Theorem]{Symmetric Square Large Sieve for  $\mathrm{PSL}_2 (\mathbb{Z} {[i]})   \backslash \mathrm{PSL}_2 (\mathbb{C}) $ \\  and Prime Geodesic Theorem for $ \mathrm{PSL}_2 (\mathbb{Z} {[i]})   \backslash \mathbb{H}^3 $}

\begin{abstract}
	  In this paper, we improve the error term in the prime geodesic theorem for the Picard manifold $ \mathrm{PSL}_2 (\mathbb{Z} {[i]})   \backslash \mathbb{H}^3  $. Instead of $ \mathrm{PSL}_2 (\mathbb{Z} {[i]})   \backslash \mathbb{H}^3  $, we establish a spectral large sieve inequality for symmetric squares over $\mathrm{PSL}_2 (\mathbb{Z} {[i]})   \backslash \mathrm{PSL}_2 (\mathbb{C}) $. This enables us to improve the bound $ O (T^{3+2/3+\vepsilon}) $ of Balkanova and Frolenkov to $ O (T^{3+1/2+\vepsilon}) $ for the second moment of symmetric square $L$-functions over $ \mathrm{PSL}_2 (\mathbb{Z} {[i]})   \backslash \mathbb{H}^3  $. The basic idea is to enlarge the spherical family $\Pi_c^{0} (T)$ of  Maass cusp forms on $ \mathrm{PSL}_2 (\mathbb{Z} {[i]})   \backslash \mathbb{H}^3 $ to the family $ \Pi_c (T, \sqrt{T}) $ of cuspidal representations on $ \mathrm{PSL}_2 (\mathbb{Z} {[i]})   \backslash \mathrm{PSL}_2 (\mathbb{C})  $. 
\end{abstract}

\author{Zhi Qi}
\address{School of Mathematical Sciences\\ Zhejiang University\\Hangzhou, 310027\\China}
\email{zhi.qi@zju.edu.cn}

\thanks{The author was supported by National Key R\&D Program of China No. 2022YFA1005300 and National Natural Science Foundation of China No. 12071420.}


\maketitle

\section{Introduction}

\subsection{Mean Lindel\"of Hypothesis and Prime Geodesic Theorem}  


Let $\varGamma = \mathrm{SL}_2 (\mathbb {Z})$ act on the hyperbolic upper half-plane $\mathbb{H}^2  $:
\begin{align*}
	\gamma z = \frac {az+b} {cz+d}, \qquad \gamma = \begin{pmatrix}
		a & b \\ c & d
	\end{pmatrix} .  
\end{align*} Let $\pi_{\varGamma} (x) $   count  the prime geodesics on $\varGamma \backslash \mathbb{H}^2$ of norm below $ x$.  
The prime geodesic theorem states
\begin{align}\label{1eq: PGT}
	\pi_{\varGamma}  (x) = \mathrm{li} (x) + E_{\varGamma}  (x), 
\end{align}
where  $\mathrm{li} (x)$ is the logarithmic integral and $E_{\varGamma} (x)$ is the error term. The Selberg trace formula directly yields  $E_{\varGamma}  (x) = O (x^{3/4+\vepsilon})$.

In the seminal work of  Iwaniec \cite{Iwaniec-PGT} in 1984, he broke the $3/4$-barrier and  obtained $ E_{\varGamma}  (x) = O (x^{35/48+\vepsilon}) $. He was the first to bring the  mean Lindel\"of hypothesis for the symmetric square (or Rankin--Selberg) $L$-functions into the study of the prime geodesic theorem. 

More explicitly, for a Hecke--Maass cusp form $u (z)$ on   $\varGamma \backslash \mathbb{H}^2$, with Hecke eigenvalues $\lambda (n)$, its symmetric square $L$-function is defined by 
\begin{align}\label{1eq: defn L(s, Sym2)}
	L (s, \mathrm{Sym}^2 u) = \zeta (2s) \sum_{n = 1}^{\infty}  \frac {\lambda (n^2)} {n^s} , \qquad (\mathrm{Re}(s) > 1), 
\end{align}
and it is well-known by the work of Shimura \cite{Shimura-Sym2} that $ L (s, \mathrm{Sym}^2 u)  $ admits analytic continuation onto the entire $s$-plane. 

By the spectral large sieve inequality of Deshouillers and Iwaniec \cite[Theorem 2]{D-I-Kuz} (see also \cite[Theorem 1]{Iwaniec-Large-Sieve}), the mean Lindel\"of hypothesis was proved later in 1995 by Luo and Sarnak \cite{Luo-Sarnak-QE}: 
\begin{align}
	\sum_{t_j \leqslant T}  | L (s, \mathrm{Sym}^2 u_j)|^2 \Lt |s|^{5+\vepsilon} T^{2+\vepsilon} , \qquad \mathrm{Re} (s) = \frac 1 2, 
\end{align}
where $u_j $ ranges over an orthonormal basis of Hecke--Maass cusp forms on $\varGamma \backslash \mathbb{H}^2$ with Laplace eigenvalue $1/4+ t_j^2$ ($t_j > 0$). This enabled them to improve the error bound into $ O (x^{7/10+\vepsilon}) $.\footnote{Note that Iwaniec also used the large sieve (for the Rankin--Selberg $L$-functions), but he used the Phragm\'en--Lindel\"of principle instead of the approximate functional equation.} 

In 2002, by inserting Burgess' estimate into Iwaniec's arguments, Cai \cite{Cai-PGT} obtained  $E_{\varGamma}  (x) = O (x^{71/102+\vepsilon})$.

The current record  $ E_{\varGamma}  (x) = O (x^{25/36+\vepsilon}) $ is due to Soundararajan and Young \cite{Sound-Young-PGT} in 2013. Their approach emphasizes the connection between prime geodesics and quadratic Dirichlet $L$-functions, and uses the Weyl subconvexity bound for the latter by Conrey and Iwaniec \cite{CI-Cubic}. 

Next, consider the Picard group $\varGamma = \mathrm{SL}_2 (\mathbb {Z}[i])$ acting on the hyperbolic space $\mathbb{H}^3 = \allowbreak  \{ (z, r) : z \in \mathbb {C}, r \in \mathbb {R}_+ \}$:
\begin{align*}
	\gamma  (z, r) = \bigg(\frac { (az+ b) (\overline {cz+d} )   +   a \overline c r^2 } {|cz+d|^2  +  |c|^2 r^2},    \frac {r  } {{|cz+d|^2 + |c|^2 r^2}}\bigg), \quad \gamma = \begin{pmatrix}
		a \! & b \\ 
		c \! & d \\
	\end{pmatrix} .
\end{align*}
As the standard reference, the reader is referred to \cite{EGM}.

In 1983, Sarnak \cite{Sarnak-H3}  established the prime geodesic theorem\footnote{More generally, let $\mathrm {F}  $ be an imaginary quadratic field (of class number $h_{\mathrm {F}} = 1$) and $\CaloO$ be its ring of integers. Sarnak's results are indeed valid for the Bianchi group  $\varGamma = \mathrm{SL}_2 (\CaloO) $, but later authors usually worked on the Picard manifold as the (spherical) Kuznetsov formula was established only explicitly over $\mathrm {F} = \mathbb {Q} (i)$ \cite{Motohashi-Kuznetsov-Picard-0}. Given the work of Lokvenec-Guleska \cite{B-Mo2} and Remark 4.2 of Koyama \cite{Koyama-PGT-Picard}, it should not be hard to extend the setting to Bianchi manifolds.} on $\varGamma \backslash \mathbb{H}^3$:
\begin{align}\label{1eq: PGT, 2}
	\pi_{\varGamma}  (x) = \mathrm{li} (x^2) + E_{\varGamma}  (x), 
\end{align}
 and he proved the benchmark error bound $ E_{\varGamma} (x) = O (x^{5/3+\vepsilon}) $.  

In 2001, motivated by the works of Iwaniec, Luo, and Sarnak, under the symmetric square mean Lindel\"of hypothesis for $\varGamma \backslash\mathbb{H}^3$,  Koyama \cite{Koyama-PGT-Picard} proved  {conditionally} that $E_{\varGamma} (x) = O (x^{11/7+\vepsilon}) $. 

Similar to \eqref{1eq: defn L(s, Sym2)}, if $ u (z, r)$ is a Hecke--Maass form on  $\varGamma \backslash\mathbb{H}^3$, with Hecke eigenvalues $\lambda (\mathfrak{n})$, where $ \mathfrak{n} $ ranges over non-zero ideals in $\CaloO = \mathbb {Z} [i]$, then 
\begin{align}\label{1eq: defn L(s, Sym2), Picard}
	L (s, \mathrm{Sym}^2 u) = \zeta_{\mathrm{F} }  (2s) \sum_{\mathfrak{n}}   \frac {\lambda (\mathfrak{n}^2)} {\RN(\mathfrak{n})^s} , \qquad (\mathrm{Re}(s) > 1), 
\end{align}
and its analytic continuation was proved in the work of Gelbart--Jacquet \cite{GJ-GL(2)-GL(3)}; here $\zeta_{\mathrm{F} }  (s)$ is the Dedekind $\zeta$ function for the Gaussian field $\mathrm{F} = \mathbb {Q} (i)$.  The mean  Lindel\"of hypothesis is of the form: \begin{align}\label{1eq: sym2 mean Lindelof}
	\sum_{t_j \leqslant T}  | L (s, \mathrm{Sym}^2 u_j)|^2 \Lt |s|^{A} T^{3+\vepsilon} , \qquad \mathrm{Re} (s) = \frac 1 2, 
\end{align}  
where $u_j $ are Hecke--Maass cusp forms on $\varGamma \backslash \mathbb{H}^3$ with Laplace eigenvalue $1+ t_j^2$ ($t_j > 0$).

However, the unconditional advance came much later in 2019: Balkanova et al.  \cite{BCCFL-PGT-Picard-1} proved the first non-trivial bound towards \eqref{1eq: sym2 mean Lindelof}: 
\begin{align}\label{1eq: sym2 mean, 1}
	\sum_{t_j \leqslant T}  | L (s, \mathrm{Sym}^2 u_j)|^2 \Lt |s|^{A} T^{4+\vepsilon} , \qquad \mathrm{Re} (s) = \frac 1 2, 
\end{align}  for an unspecified constant $A > 0$,   and,  by the argument of Koyama, achieved $E_{\varGamma} (x) = O (x^{13/8+\vepsilon})$. Indeed, it is implicit in the work of  Koyama (see the discussion below Corollary 3.4 in \cite{BCCFL-PGT-Picard-1}) that  
\begin{align}\label{1eq: bound for E, Koyama}
E_{\varGamma} (x) \Lt x^{\frac {11+4\phi} {7+2\phi} + \vepsilon  } ,
\end{align}  if one had 
\begin{align}\label{1eq: mean bound}
	\sum_{t_j \leqslant T}  | L (s, \mathrm{Sym}^2 u_j)|^2 \Lt |s|^{A} T^{3+2 \phi +\vepsilon} , \qquad \mathrm{Re} (s) = \frac 1 2.  
\end{align} 

More recent improvements in 2022 rely on the subconvexity of quadratic Dirichlet $L$-functions over $\CaloO$, that is,
\begin{align}\label{1eq: subconvex bound}
	L (s, \vchi_{\mathfrak{d}} ) \Lt |s|^A \RN (\mathfrak{d})^{\theta +\vepsilon}, \qquad \mathrm{Re} (s) = \frac 1 2, 
\end{align}
where $ \vchi_{\mathfrak{d}} $ is the quadratic Dirichlet character of conductor $\mathfrak{d}$; the currently best Weyl subconvexity exponent $\theta = 1/6$ is due to Nelson \cite{Nelson-Eisenstein}. In the spirit of Soundararajan and Young \cite{Sound-Young-PGT} (see also \cite{C-Wu-Z-KB-Formula}), it is proved by Balog et al. \cite{BBCL-PGT-Picard-3} (stated more explicitly in \cite[\S 6]{BF-Picard-Sym2}) that if both \eqref{1eq: mean bound} and \eqref{1eq: subconvex bound} hold then \eqref{1eq: bound for E, Koyama} may be improved into
\begin{align}\label{1eq: error bound}
	E_{\varGamma} (x) \Lt x^{\frac 3 2 + \frac {24 \theta -1 +  (16\theta+2) \phi} {46+20\phi} + \vepsilon  } ,
\end{align}
and hence $ E_{\varGamma} (x) = O (x^{3/2+4\theta/7+\vepsilon}) $ on $\phi = 1/2$ and $O (x^{67/42+\vepsilon})$ on $\theta = 1/6$.\footnote{Note that $E_{\varGamma} (x) = O (x^{3/2+\theta+ \vepsilon})$ was obtained in \cite{BF-PGT-Picard-2}, and it recovers Sarnak's bound on the Weyl exponent $\theta=1/6$.}
Then Balkanova and Frolenkov  \cite{BF-Picard-Sym2} successfully proved
\begin{align}\label{1eq: sym2 mean}
	\sum_{t_j \leqslant T}  | L (s, \mathrm{Sym}^2 u_j)|^2 \Lt |s|^{A} T^{3+ 4 \theta +\vepsilon} , \qquad \mathrm{Re} (s) = \frac 1 2, 
\end{align}  
and the currently best bound  
\begin{align}
	 E_{\varGamma} (x) \Lt x^{3/2+41/474+\vepsilon}    , 
\end{align} 
on $\phi = 1/3$ and $\theta = 1/6$.

In this paper, we improve the results of Balkanova and Frolenkov into
\begin{align}\label{1eq: mean bound, 1}
	\sum_{t_j \leqslant T}  | L (s, \mathrm{Sym}^2 u_j)|^2 \Lt |s|^{3+\vepsilon} T^{3 +1/2 +\vepsilon} , \qquad \mathrm{Re} (s) = \frac 1 2, 
\end{align}   
which goes halfway between \eqref{1eq: sym2 mean, 1} and
the mean-Lindel\"of bound, 
and
\begin{align}
	E_{\varGamma} (x) \Lt x^{3/2+25/306+\vepsilon}  ,
\end{align} 
on $\phi = 1/4$ and $\theta = 1/6$. Note that $3/2+ 41/474 = 1.58649... $  while $ 3/2+25/306 = 1.58169... $. 
 
These results will be proved via the original large sieve approach by Iwaniec, Luo, and Sarnak  \cite{Iwaniec-PGT,Luo-Sarnak-QE}. However, it is remarked in \cite{BCCFL-PGT-Picard-1}:   
		``In contrast to the proof in Luo and Sarnak, applying the best known large sieve in $\mathbb{H}^3$ (of Watt) is not helpful.'' 
Our idea to resurrect the large sieve approach is to rise from the ground $\mathbb{H}^3 \cong \allowbreak \mathrm {SL}_2 (\mathbb {C}) /  \mathrm{SU}_2 (\mathbb {C})$ and go into the space $\mathrm {SL}_2 (\mathbb {C})$.

The spectral large sieve for $\mathrm {SL}_2 (\mathbb {Z}[i])\backslash \mathrm {SL}_2 (\mathbb {C})$ of Watt \cite{Watt-1-Large-Sieve} was recently improved by the author \cite{Qi-Spectral-LS} so that it now has the strength of that of Deshouillers and Iwaniec \cite{D-I-Kuz} in the most `balanced' case (the spherical case is the most `unbalanced' case, so the remark of Balkanova et al. is alright---the large sieve restricted to $ \mathbb{H}^3 $ is not strong enough). It will be explained in a moment that this already enables us to recover \eqref{1eq: sym2 mean, 1}.   

The  spectral large sieve considered in this paper is for the symmetric squares on $\mathrm {SL}_2 (\mathbb {Z}[i])\backslash \mathrm {SL}_2 (\mathbb {C})$. It will be established on the analysis of Bessel integral in  \cite{Qi-Spectral-LS} and the evaluation of certain character sums in Young's recent work \cite{Young-Sym2-LS}. 


\subsection{Main Results} \label{sec: main results}


Let $\mathrm {F} = \mathbb {Q} (i)$ and $ \CaloO  = \mathbb {Z} [i]$.  Subsequently the letter $n$ will always stand for a non-zero Gaussian integer  
 in $\CaloO    $. Let 
 $ (n)$ denote the ideal $n \CaloO$. Let $\RN (n) = |n|^2  = |\CaloO / n \CaloO|$ denote the norm of $n$ or $(n)$. 

Let $\varGamma = \mathrm{PSL}_2 (\CaloO)$ and $G = \mathrm{PSL}_2 (\mathbb {C})$. In this paper, we shall be mainly concerned with the discrete spectrum $\Pi_{c} $ of  the irreducible cuspidal representations  in $L^2  (\varGamma \backslash G)$. 

For $\uppii \in \Pi_c$, let $  (t      _{\uppii}, p          _{\uppii}) $ be the spectral parameter of $\uppii$ in the sense that $\uppii$ is isomorphic to the principal series of $G$ unitarily induced from the character
\begin{align*}
   \begin{pmatrix}
		\sqrt{z} & v \\
		& 1/\sqrt{z}
	\end{pmatrix}   \rightarrow 
	|z|^{i t      _{\uppii} } (z/|z|)^{  p          _{\uppii}}, 
\end{align*} where $t_{\uppii}$ is real  and $p          _{\uppii}$ is integral. Let $\lambda_{\uppii} (n)$ denote the Hecke eigenvalues of $ \uppii $. We have $  \lambda_{\uppii} (-n) = \lambda_{\uppii} (n) $ so that $\lambda_{\uppii} (n^2)$  only depends on the ideal $  (n)$. 
Let $L (s, \mathrm{Sym}^2 \uppii)$ be the   symmetric square $L$-function attached to $\uppii$. Define the harmonic weight
\begin{align*}
	\omega_{\uppii} = \frac 1 {L (1, \mathrm{Sym}^2 \uppii)}.  
\end{align*}

Define
\begin{align*}
	  \Pi_{c} (T, P) =  \big\{ \uppii \in \Pi_{c} : |t      _{\uppii} | \leqslant T,   |p          _{\uppii} | \leqslant P \big\},
\end{align*}
and for a sequence $ \boldsymbol{a}: \CaloO / \CaloO^{\times} \rightarrow \mathbb {C} $ (the quotient $\CaloO / \CaloO^{\times}$ is identified with the set of ideals)
\begin{align}\label{1eq: defn of C}
	\SC^{(2)} (T, P, N) =	 \sum_{\uppii \shskip \in \Pi_{c} (T, P) }  \omega_{\uppii} \bigg| \sum_{ \RN (n)  \leqslant N }  a_{(n)} \lambda_{\uppii} (n^2)  \bigg|^2. 
\end{align}

The main result of the recent paper  \cite{Qi-Spectral-LS} may be modified\footnote{The spectral large sieve inequality in  \cite{Qi-Spectral-LS} is on $ \mathrm{PGL}_2 (\CaloO) \backslash \mathrm{PGL}_2 (\mathbb {C}) $ so that $ \lambda_{\uppii} (n) $ is  dependent only on the ideal $(n)$ as required by the definition of $L (s, \uppii)$. However,  the proof therein may be easily adapted to the setting of $ \mathrm{PSL}_2 (\CaloO) \backslash \mathrm{PSL}_2 (\mathbb {C}) $. } as the following spectral large  sieve inequality: 
\begin{equation}\label{1eq: large sieve, 0}
	\begin{split}
		\sum_{\uppii \shskip \in \Pi_{c} (T, P) }  \omega_{\uppii}   \bigg|   \sum_{  \RN(n)  \leqslant   N }   a_{n} \lambda_{\uppii} (n)  \bigg|^2   \Lt    \big(T^2 + P^2\big) \bigg( T P +   \frac{ N } {T P}  \bigg) (T PN)^{\vepsilon}   \sum_{  \RN(n)  \leqslant  N } |a_{n}|^2,
	\end{split}
\end{equation}
for any sequence $ \boldsymbol{a}: \CaloO / \{\pm 1 \} \rightarrow \mathbb {C} $. This improves the bound of Watt \cite[Theorem 1]{Watt-1-Large-Sieve} where it is $N / \sqrt{T P}$ instead of $N / T P$. A direct consequence of  \eqref{1eq: large sieve, 0} is the bound 
\begin{equation}\label{1eq: sym2 large sieve, 0}
	\begin{split}
	\SC^{(2)} (T, P, N)    \Lt    \big(T^2 + P^2\big) \bigg( T P +   \frac{ N^2 } {T P}  \bigg) (T PN)^{\vepsilon}   \sum_{  \RN(n)  \leqslant  N } |a_{(n)}|^2,  
	\end{split}
\end{equation} 
and in the  square case $P = T$ it reads
\begin{equation}\label{1eq: sym2 large sieve, 0.1}
	\begin{split}
		\sum_{\uppii \shskip \in \Pi_{c} (T, T) }  \omega_{\uppii} \bigg| \sum_{ \RN (n)  \leqslant N }  a_{(n)} \lambda_{\uppii} (n^2)  \bigg|^2  \Lt     \big( T^4 +    { N^2 }    \big) (TN)^{\vepsilon}   \sum_{  \RN(n)  \leqslant  N } |a_{(n)}|^2.
	\end{split}
\end{equation}  

Now let $\Pi_{c}^0$ denote the set of spherical $\uppii$ with $p          _{\uppii} = 0$, and, by convention,  identify $\Pi_{c}^0$ with the basis of Maass cusp forms $u_j $ as above. Define
\begin{align*}
	\Pi_{c}^{0} (T ) = \big\{ u_j \in \Pi_{c}^0 : t      _{j}   \leqslant T  \big\}. 
\end{align*}

Our observation is that if we embed $ \Pi_c^0 (T) \hookrightarrow \Pi_c (T, T)$ and apply \eqref{1eq: sym2 large sieve, 0.1} with $N \asymp_{s} T^{2+\vepsilon}$, then the moment bound $O_s (T^{4+\vepsilon})$ as in \eqref{1eq: sym2 mean, 1}  is  readily recovered; it is indeed the mean-Lindel\"of bound for the larger family $\Pi_c (T, T)$! 
It turns out that the benchmark bound \eqref{1eq: sym2 large sieve, 0} may be improved by exploiting the arithmetical feature of the inner $n$-sum in \eqref{1eq: defn of C}, namely the square $n^2$ in  $\lambda_{\uppii} (n^2)$. Consequently, we can verify the mean-Lindel\"of bound for the family  $\Pi_c (T, \sqrt{T})$ and thereby prove the moment bound $O_{s} (T^{3+1/2+\vepsilon})$ as in \eqref{1eq: mean bound, 1} for the spherical family $\Pi_c^{0} (T)$ via the embedding $ \Pi_c^0 (T) \hookrightarrow \Pi_c (T, \sqrt{T}) $. 

\begin{thm}\label{thm: sym2 large sieve}
	 Let $T, P , N  \geqslant  1/2$ {\rm(}say{\rm)} and $\vepsilon > 0$ be real numbers.  Let $\boldsymbol{a} : \CaloO / \CaloO^{\times}   \rightarrow \mathbb {C}$ be a sequence. Then 
	 \begin{equation}\label{1eq: sym2 large sieve, 1} 
	 	\begin{split}
	 			\SC^{(2)} (T, P, N) \Lt  \bigg\{   T & P   \big(T^2 + P^2 \big) + T P N + \\
	 		& + \frac {T^2+P^2} {T P} \bigg(\frac 1 {T^2} + \frac 1 {P^2} \bigg) N^2 \bigg\} (T PN)^{\vepsilon} \sum_{  \RN(n)  \leqslant  N } |a_{(n)}|^2. 
	 	\end{split}
	 \end{equation}
\end{thm}

Similar to \cite[Theorem 1]{Qi-Spectral-LS}, the estimate in \eqref{1eq: sym2 large sieve, 1}  also applies to the contribution from the Eisenstein series (see \eqref{2eq: defn of E}), but here it is not our object of concern.  

The reader may compare \eqref{1eq: sym2 large sieve, 0.1} with  
 \begin{equation} 
 	\sum_{\uppii \shskip \in \Pi_{c} (T, \sqrt{T}) }  \omega_{\uppii} \bigg| \sum_{ \RN (n)  \leqslant N }  a_{(n)} \lambda_{\uppii} (n^2)  \bigg|^2 \Lt    \bigg(   T^{7/2} +     \frac {N^2} {T^{1/2}}   \bigg) (TN)^{\vepsilon} \sum_{  \RN(n)  \leqslant  N } |a_{(n)}|^2, 
 \end{equation}
which is a rewriting of \eqref{1eq: sym2 large sieve, 1} in the case $P = \sqrt{T}$.

%

\begin{thm}\label{thm: mean-Lindelof}
	Let $\mathrm{Re} (s) = 1/2$. We have
	\begin{equation}
		\sum_{\uppii \shskip \in \Pi_{c} (T, P) }   |L(s , \mathrm{Sym}^2 \uppii )|^2 \Lt |s|^{3+\vepsilon} \cdot \left\{ 
		\begin{split}
			& (P + \sqrt{T}) T^{3+\vepsilon} , & & \text{ if }  P \leqslant T, \\
			&  (T+\sqrt{P}) P^{3+\vepsilon}, & & \text{ if }  T \leqslant P,
		\end{split} \right. 
	\end{equation}
and hence 
\begin{align}\label{1eq: mean sym2}
	\sum_{t_j \leqslant T }   |L(s , \mathrm{Sym}^2 u_j )|^2 \Lt |s|^{ 3+\vepsilon  } T^{3+1/2+\vepsilon},   
\end{align}
 for any $\vepsilon > 0$. 
\end{thm}

The next theorem follows directly from \eqref{1eq: PGT, 2}, \eqref{1eq: mean bound}--\eqref{1eq: error bound}, and \eqref{1eq: mean sym2}.  
\begin{thm}
We have 
	 \begin{align}
	 	   	\pi_{\varGamma}  (x) = \mathrm{li} (x^2) + O \big( x^{ 3/2 + {(56 \theta - 1)} / {102}    + \vepsilon } \big) , 
	 \end{align} 
 for any $\vepsilon > 0$. 
\end{thm}



For technical simplifications, 
in the main body of this paper, we shall replace the range $ \RN (n) \leqslant N $ as in the statements above  by the dyadic  $ N < \RN (n) \leqslant 2 N$ (occasionally written as $ \RN (n) \sim N $). It is clear that the results so modified are of the same strength. Subsequently, let
\begin{align*}
		\| \boldsymbol{a}_N \|_2   = \bigg(  \sum_{  N <  \RN(n)  \leqslant   2 N } |a_{(n)}|^2 \bigg)^{1/2}. 
\end{align*}

\subsection{Quadratic Forms with Kloosterman Sums} 

Let $e [z] = \exp (2\pi i\mathrm{Re} (z))$. 
For $  m,  n \in \CaloO $ and $c \in \CaloO \smallsetminus \{0\}$,      define the Kloosterman sum
\begin{align}\label{1eq: defn Kloosterman}
	S   (m, n ; c ) = \sumx_{   \valpha  (\mathrm{mod} \, c) } e \bigg[  \frac {  \valpha m +   \bar{\valpha} n } {c} \bigg] ,
\end{align}  
where the {\scriptsize \ding{72}} indicates the condition $(\valpha, c) = (1)$ and $ \bar{\valpha}$ is given by $\valpha \bar{\valpha}  \equiv 1 (\mathrm{mod} \, c)$. For  $(w, c) = (1)$, define
\begin{align}\label{1eq: defn of F}
	F  (w; c) = S (w^2, 1; c) e  [2 w / c ]. 
\end{align}
The $e  [2 w / c ]$ will be extracted from the Bessel integral---it is an important feature in the symmetric-square problems. 
 
	\begin{prop}\label{prop: quad form}
		
		Let $\boldsymbol{a}, \boldsymbol{b} : \CaloO / \CaloO^{\times}   \rightarrow \mathbb {C}$ be a pair of sequences,  $ d \in \CaloO \smallsetminus \{0\}$, and $\theta \in \mathbb {C} \smallsetminus \{0\} $. Let $C, M, N \geqslant 1/2$, $\vepsilon > 0$, and $\gamma$ be real. Set  \begin{equation*}
			K = C + \sqrt{CMN}|\theta|. 
		\end{equation*}   
		Write
		\begin{align*}
			F^{\gamma} ( d, \theta,   C, M, N) = \mathop{\mathop{\sum \sum}_{ {M} < \RN (m) \leqslant   2 M } }_{N < \RN (n) \leqslant 2 N}  a_{(m)}  \overline{b }_{(n)}  \mathop{\sum_{C < \RN (c) \leqslant 2 C}}_{(c, d mn) = (1)}   {\RN(c)^{\gamma}}  {F  (d m n; c) e [  mn \theta / c]} .
		\end{align*} 
	Then 
		\begin{equation} 
			\label{1eq: Kloosterman} 	 \begin{split}
				  F^{\gamma} ( d, \theta,   C, M, N)   \Lt_{\vepsilon}  C^{1+\gamma} \big( K + \sqrt{M } + \sqrt{ N} + C \sqrt{MN}  /K\big)    {K^{\vepsilon}}   \|\boldsymbol{a}_{M} \|_2  \|\boldsymbol{b}_{N} \|_2 &.
			\end{split}
		\end{equation} 
\end{prop}

Proposition \ref{prop: quad form} is in a similar spirit to Proposition 3 in \cite{D-I-Kuz} and Proposition 2 in \cite{Qi-Spectral-LS}. However, here we also need to separate $m$ and $n$ in $F (dmn; c)$ by the non-Archimedean Mellin technique. Fortunately, this has been done (in the $\mathrm {SL}_2 (\mathbb {Z})$ setting) recently by Young \cite{Young-Sym2-LS}.

\subsection{Remarks}\label{sec: remarks}


The unpolished symmetric-square large sieve for  $\mathrm {SL}_2 (\mathbb {Z})$ of Young in  \cite[(1.32)]{Young-Sym2-LS} reads
\begin{align}
	\sum_{T \leqslant t_j \leqslant T+\varDelta} \omega_j \bigg| \sum_{ n \leqslant N }  a_{n} \lambda_{j} (n^2)  \bigg|^2 \Lt   \bigg(\varDelta T + \varDelta N + \frac {N^2} {\varDelta} \bigg) (T N)^{\vepsilon} \sum_{ n \leqslant N }  |a_{n}|^2,  
\end{align}
for $\varDelta \leqslant T$, while our \eqref{1eq: sym2 large sieve, 1} reads
\begin{align}\label{1eq: sym2 P<T}
\sum_{\uppii \shskip \in \Pi_{c} (T, P) }   \omega_{\uppii}   \bigg|   \sum_{  \RN(n)  \leqslant   N }    a_{(n)} \lambda_{\uppii} (n^2)  \bigg|^2  \Lt 	\bigg( T^3P + T P N +   \frac {T N^2} {P^3} \bigg) (TN)^{\vepsilon}   \sum_{  \RN(n)  \leqslant  N }  |a_{(n)}|^2, 
\end{align} 
for $P \leqslant T$. The reader may see the clear similarity, but it is more interesting to note that $\sqrt{T}$ is the barrier for both to achieve the mean Lindel\"of bound for the symmetric square $L$-functions. {However,  some major differences would be found if one analyzes the Bessel integrals closely. The  family truly analogous to 
	\begin{align*}
		\big\{  u_j : T \leqslant t_j \leqslant T+ \varDelta  \big\} 
	\end{align*}
	should be
	\begin{align*}
		\big\{ \uppii :  T \leqslant |t_{\uppii}|, \, |p_{\uppii}| \leqslant T+ \varDelta   \big\}. 
\end{align*}}

For the moment problem of $L$-functions, it is understood that the approach via large sieve can usually be surpassed by that via Poisson/Vorono\"i summation. 
In another paper of Young with Khan \cite{Khan-Young-Sym2}, the mean Lindel\"of bound is proved for $\varDelta = T^{1/5 }$ by Poisson summation (thrice) and the quadratic large sieve of Heath-Brown \cite{Heath-Brown-Quad-LS}. It seems that in our setting their approach may be used to settle the mean Lindel\"of  hypothesis for all $P$ and $ T$. 

\subsection{Refinement}

In the same fashion as \cite[Theorem 1.2]{Young-Sym2-LS}, we may improve \eqref{1eq: sym2 large sieve, 1} into
\begin{align} 
	\SC^{(2)} (T, P, N) \Lt  \big(T^2+P^2\big) \big(TP + {\textstyle \sqrt[4]{T^2+P^2 }}\sqrt{N} \big) (TP)^{\vepsilon} \sum_{  \RN(n)  \leqslant  N } |a_{(n)}|^2  ,
\end{align}
if $N \leqslant T^2+ P^2$, 
\begin{align} 
	\SC^{(2)} (T, P, N) \Lt  \big(TP + {\textstyle \sqrt{T^2+P^2}} \sqrt[4]{N} \big) N (TP)^{\vepsilon} \sum_{  \RN(n)  \leqslant  N } |a_{(n)}|^2  ,
\end{align}
if $T^2 + P^2 < N \leqslant T^4+P^4$, and 
\begin{align} 
	\SC^{(2)} (T, P, N) \Lt \frac {N^2} {T^2+P^2} N^{\vepsilon} \sum_{  \RN(n)  \leqslant  N } |a_{(n)}|^2  ,
\end{align}
if $N > T^4 + P^4$. 

Say $ P \leqslant T$ so that \eqref{1eq: sym2 large sieve, 1} is simplified into \eqref{1eq: sym2 P<T}. Since $\SC^{(2)} (T, P, N)$ on the left is non-decreasing in $P$, these improvements may be seen by enlarging $P$ as follows: for $N \leqslant T^2$, replace $P$ by $P+\sqrt{N/T}$; for  $T^2 < N \leqslant T^4$,  replace $P$ by $P+\sqrt[4]{N}$; for $N > T^4$,  replace $P$ by $P+T$.  

 \subsection*{Convention}
 
Throughout the paper,  $\vepsilon  $ is arbitrarily small and its value  may differ from one occurrence to another.


\noindent \subsection*{Acknowledgments} The author thanks the referees for their careful readings and helpful comments. 



\section{Basic Notions}\label{sec: notation} 

Denote $\mathrm {F} = \mathbb {Q} (i)$ and $ \CaloO  = \mathbb {Z} [i]$. Let $\CaloO^{\times}$ be the group of units and $\CaloO    / \CaloO^{\times}$ denote the set of   ideals. 
Reserve  $m$, $n$, and $c$  for   non-zero  Gaussian integers. 
Use the corresponding Gothic letters  $\frm = (m)$,  $\frn = (n)$,  $\frc = (c)$...... for non-zero ideals in $\CaloO$. Let  $\RN (n) = |n|^2$ or $\RN (\frn) = |\CaloO / \frn| $ denote the norm of $n$ or $\frn$.  

As usual, let $\tau (\frn)$, $\varphi (\frn)   $, and $\mu (\frn)$  respectively  denote the divisor function,  Euler totient function, and M\"obius function over $\CaloO$.


For a  function $F  $ defined on $(\CaloO / \frc)^{\times}$,  we have the Mellin inversion formula:
\begin{align}\label{2eq: Mellin inversion}
	F (w) = \sum_{\vchi (\mathrm{mod}\, \frc)}  \widehat{F} (\vchi) \vchi (w), \qquad    \widehat{F} (\vchi) = \frac 1 {\varphi (\frc) } \sumx_{\valpha (\mathrm{mod}\, \frc)} \overline{\vchi} (\valpha) F (\valpha), 
\end{align}
where the sum of $\vchi$ is over $  \widehat {(\CaloO / \frc)}{}^{\times}$, namely, $\vchi$ are characters  of the group  $(\CaloO / \frc)^{\times}$.  It will be easier to consider $\vchi (n)$ as a Dirichlet character on $ \CaloO $ of modulus $\frc$.  

The   (unitary) characters on $\mathbb {C}^{\times} $ are of the form     \begin{align}
	\vlambda_{i t      ,  p          } (z) = |z|^{  i t        } (z/|z|)^{  p          }  , 
\end{align}  for   $ t         $ real and   $ p           $ integral, 
so we introduce  $\boldsymbol{\EuScript{A}}     = \mathbb {R} \times  \mathbb {Z}$ to parameterize the (Mellin) unitary dual of $\mathbb {C}^{\times}  $.  
Let $\mathrm{d}\mu (t      ,  p) $ 
denote the usual Lebesgue measure on  $\boldsymbol{\EuScript{A}}    $. 
For simplicity, write $ 	\vlambda_{  p          } (z) = 	\vlambda_{0,  p          } (z)$.

Moreover, let
\begin{align*}
	\boldsymbol{\EuScript{A}} (T, P) =   \big\{ (t      , p          ) \in \boldsymbol{\EuScript{A}} : |t        | \leqslant T,   |p            | \leqslant P \big\}. 
\end{align*}


For $p           \in \mathbb {Z} $,  define 
\begin{align}\label{1eq: zeta}
	\zeta (s, p          ) =  \sum_{\frn \shskip \subset \CaloO } \frac {	\vlambda_{4 p          } (\frn) } {\RN(\frn)^{s}}, \qquad \text{($\mathrm{Re} (s) > 1$)}, 
\end{align}
to be the Hecke $\zeta$ function associated to the Gr\"ossencharakter $	\vlambda_{4 p          } : (n) \rightarrow (n/|n|)^{4 p          }$, and   define 
\begin{align}\label{1eq: tau s}
	\tau_{s, p          } (\frn) = \sum_{ \mathfrak{a} \mathfrak{b} = \frn} 	\vlambda_{4 p          } (\mathfrak{a} \mathfrak{b}^{-1}) \RN ( \mathfrak{a} \mathfrak{b}^{-1} )^{s} .
\end{align}  
It is well known that $ \zeta (s, p          ) $ has analytic continuation to the whole complex plane except for a simple pole at $s = 1$ in the case $p           = 0$ as it is the Dedekind $\zeta$ function $\zeta_{\mathrm {F}} (s)$.  These definitions will only be used to make explicit the Eisenstein contribution. For example, the Eisenstein counterpart of $ 	\SC^{(2)} (T, P, N) $ as in \eqref{1eq: defn of C} should read:
\begin{align}\label{2eq: defn of E}
	\SE^{(2)} (T, P, N) =	\int \!\!\! \int_{ \boldsymbol{\EuScript{A}} (T/2, P/4) }  \omega (t, p          ) \bigg|  \sum_{ \RN (n)  \leqslant N }   a_{(n)} \tau_{i t      , p          } (n^2)  \bigg|^2  \mathrm{d}\mu (t      , p           ),
\end{align}
with
\begin{align*}
\omega (t, p          ) =	\frac 1 { | \zeta (1+2 i t      , 2 p          )   |^2  }. 
\end{align*}


As before, let $\varGamma = \mathrm{PSL}_2 (\CaloO)$ and $G = \mathrm{PSL}_2 (\mathbb {C})$.

For $(t      , p          ) \in \boldsymbol{\EuScript{A}}$, let $\uppii_{i t      , p          } $ be the unitary principal series of $G$: the unique infinite dimensional constituent of the representation   unitarily induced from the character
\begin{align*}
	\vlambda_{i t      ,  p          }: 	\begin{pmatrix}
		\sqrt{z} & v \\
		& 1/\sqrt{z}
	\end{pmatrix}   \rightarrow 
	|z|^{i t      } (z/|z|)^{  p          }. 
\end{align*}
Note that  $\uppii_{-i t      , - p          } \cong   \uppii_{i t      , p          } $.
As the Selberg conjecture holds for $\varGamma$, we do not consider complementary series (\cite[\S 6.7.6]{EGM}).  

Let  $  L^2_{c} (\varGamma \backslash G )  $ be the space of  cusp forms for $ \varGamma $.  Let $\Pi_{c} $ denote the discrete   spectrum of the irreducible constituents of $L^2_{c} (\varGamma \backslash G )$. 

It may be assumed that each  $\uppii \in \Pi_{c}$ is Hecke invariant. Let $\lambda_{\uppii} (n)$ be the Hecke eigenvalues of $ \uppii $. It is known that $\lambda_{\uppii} (n)$ are all real.  Note that the $\varGamma$-invariance implies that $ \lambda_{\uppii} (- n) = \lambda_{\uppii} (n) $ (consider the diagonal matrices in $\varGamma = \mathrm{PSL}_2(\CaloO) $).

For $\uppii \in \Pi_{c}$, let $  (t      _{\uppii}, p          _{\uppii})$ be the parameter $(t      , p          ) \in \boldsymbol{\EuScript{A}}$ such that  $ \uppii \simeq \uppii_{i t       , p           }$.

Define 
\begin{align*}
	\Pi_{c} (T, P) =  \big\{ \uppii \in \Pi_{c} : |t      _{\uppii} | \leqslant T,   |p          _{\uppii} | \leqslant P \big\} . 
\end{align*}

\section{Gallagher's Hybrid Large Sieve}

The hybrid large sieve of Gallagher \cite[Theorem 2]{Gallagher-LS} has been extended to number fields by Duke \cite[Theorem 1.1]{Duke-LargeSieve}. 

\begin{lem}\label{lem: hybrid}
Let $C, T, N  \geqslant  1 $  and $\vepsilon > 0$ be real numbers.  Let $\boldsymbol{a} : \CaloO / \CaloO^{\times}   \rightarrow \mathbb {C}$ be a sequence. Then 
\begin{equation} \label{2eq: hybrid} 
\begin{split}
	\sum_{\RN(\frc) \leqslant C}  \, \sumst_{\vchi (\mathrm{mod}\, \frc)}  	\int \!\!\! \int_{ \boldsymbol{\EuScript{A}} (T, T) }   \bigg|   \sum_{ \RN(n) \leqslant N }  a_{(n)} & \vchi (n)\vlambda_{it      ,  p          } (n)  \bigg|^2 \mathrm{d}\mu (t      , p)\\
	& \Lt  \big( C^2 T^2 + N \big) (CT)^{\vepsilon} \sum_{ \RN(n) \leqslant N }  |a_{(n)}|^2 ,
\end{split}
\end{equation}
where the {\scriptsize \ding{72}}   restricts to primitive $\vchi$.  
  
\end{lem}

	Actually, Duke's large sieve is for Gr\"ossencharaktere (instead of Dirichlet characters), so one needs to impose the condition $ \vchi (\epsilon) \lambda_{p          } (\epsilon) = 1  $ ($\epsilon \in \CaloO^{\times}$), but this may be safely removed as the $n$-sum vanishes if otherwise. 

\section{Proof of Proposition \ref{prop: quad form}}

\subsection{Non-Archimedean Analysis of Young} 

For brevity, let us write $F (w) = F(w; c) $. First, we recollect some results of Young  in \cite[\S 4]{Young-Sym2-LS} on the character sum \begin{align}\label{4eq: F (chi)}
	\widehat{F}  (\vchi) = \frac 1 {\varphi (c) } \, \sumst_{\valpha (\mathrm{mod}\, c)} \bar{\vchi} (\valpha) F  (\valpha), \qquad F  (w) =  S ( w^2, 1; c) e \bigg[\frac {2  w} c \bigg]; 
\end{align} 
his local calculations may be readily extended to any non-Archimedean field.   

In this  section, for a Dirichlet character $\vchi$ we always let $\mathfrak{c} $ be its modulus and $\mathfrak{c}^{\star}  $ be its conductor, and later if $\vchi$ had a subscript then it will also be added to $\frc$. 

Let $\vchi_0$ denote the trivial character. By definition $\vchi$ is primitive if $\mathfrak{c} = \mathfrak{c}^{\star} $. We say $\vchi$ is {semi-primitive} if for every prime $ \frp | \frc$ we have $1 \leqslant v_{\frp} (\frc^{\star}) < v_{\frp} (\frc)$, where $v_{\frp}$ is the $\frp$-adic valuation.

{Note that $\widehat{F}  (\vchi)$  depends mildly on the choice of $c$ up to $\CaloO^{\times}  $ (as $F ( w ; \epsilon c) = F (\bar{\epsilon} w;  c)$), but this will not affect the analysis.}

By the Chinese remainder theorem, $\widehat{F} (\vchi)$ satisfies the twisted-multiplicative relation:
\begin{align}\label{4eq: mult}
	\widehat{F}(\vchi_1 \vchi_2 ) = \overline{\vchi}_1 (c_2)  \overline{\vchi}_2 (c_1) \widehat{F}(\vchi_1  ) \widehat{F}(  \vchi_2 ), \qquad (c_1, c_2) = (1), 
\end{align}
so $  | 	\widehat{F} (\vchi)  |$ is multiplicative.   Thus it suffices to consider the prime-power case $\frc = \frp^k$.

\begin{lem}\label{lem: Mellin F (chi)}
	Let the notation be as above.

	{\rm(1. Primitive Case)}  If $  \frc = \frc^{\star} = \frp^k$ for $k \geqslant 1$, then  for $ \RN(\frp) = 2 $ {\rm(}$\frp = \allowbreak(1+i)${\rm)} $ 	\widehat{F} (\vchi)  = 0$, while, for $ \RN(\frp) $ odd, 
	\begin{equation}\label{4eq: primitive}
\big| 	\widehat{F} (\vchi) \big|=	\left\{  	\begin{split} 
&	\displaystyle \frac {\sqrt{\RN (\frp)}}   {\RN(\frp) - 1},    & &   \vchi^2 = \vchi_0, 
\\
	&	\displaystyle \frac {\RN (\frp)} {\RN(\frp) - 1},    & &   \vchi^2 \neq \vchi_0 .
		\end{split} \right. 
	\end{equation} 

{\rm (2. Trivial Case)} If $\vchi = \vchi_0 $ is trivial of modulus $\frp^{k}$, then
\begin{equation}\label{4eq: trivial}
	 	\widehat{F} (\vchi_0)  =	\left\{  	\begin{split} 
		&	\displaystyle \frac {1} {\RN(\frp) - 1},    & &   k = 1,   \\
		&	\displaystyle  \RN(\frp)^{k/2},    & &   k \text{ even},  \\
		& 0, & &     \text{else}. 
	\end{split} \right. 
\end{equation}
	 
{\rm(3. Semi-Primitive Case)} 	 Suppose that $\frc = \frp^{k}$ and $\frc^{\star} = \frp^{k^{\star}}$ with $ 1 \leqslant k^{\star} < k $. If  $k^{\star}  \not\equiv k (\mathrm{mod}\, 2) $ then $ \widehat{F} (\vchi) = 0 $. If $ k^{\star}  \equiv k (\mathrm{mod}\, 2) $, then 
\begin{equation}\label{4eq: semi-primitive, 1}
	\big|\widehat{F} (\vchi)\big| \leqslant \RN(\frp)^{k/2}, 
\end{equation}
for $\vchi^2  = \vchi_0$, and
\begin{equation}\label{4eq: semi-primitive, 2}
	\big|\widehat{F} (\vchi) \big|  =	\left\{  	\begin{split} 
		&	 O(  1),    & &    \RN (\frp) = 2,     \\
		&	0 ,    & &     \RN (\frp) \text{  odd}, 
	\end{split} \right. 
\end{equation}
for $\vchi^2  \neq \vchi_0$. 

\end{lem}

\begin{remark}
Note that one necessarily has $k = 1$ or $k^{\star} = 1$ 
in the first  case of \eqref{4eq: primitive} or  \eqref{4eq: semi-primitive, 1} when $\vchi^2  = \vchi_0$ respectively. 
\end{remark}

\begin{remark}
	 Some remarks are in order on {\rm(}the local calculations for{\rm)} Lemma \ref{lem: Mellin F (chi)}. Similar to {\rm(3.8)} in \cite{Young-Sym2-LS}, the basic formula for $  \widehat{F}  (\vchi) $ reads{\rm:}
	 \begin{align*}
	 	\widehat{F}  (\vchi) = \frac 1 {\varphi (c) } \, \mathop{\sumst \sumst}_{\valpha, \beta (\mathrm{mod}\, c)} \overline{\vchi} (\valpha \beta) e   \bigg[\frac {(\valpha+1)^2 \beta }  c \bigg]. 
	 \end{align*}
	 
	
	For the primitive case, the $\beta$- and the $\valpha$-sums are evaluated by Gauss and Jacobi sums, so the proof of Lemma 4.2 in \cite{Young-Sym2-LS} may be carried out here if the formula {\rm(3.18)} in \cite{IK} for Gauss and Jacobi sums were established on $\CaloO / \mathfrak{c} $. 
	
	For the trivial case, the $\beta$-sum turns into a Ramanujan sum, and the proof of Lemma 4.3 in \cite{Young-Sym2-LS} may be applied here almost identically. 
	
	For the semi-primitive case, the situation is more delicate, but the proof of  Lemma 4.6 in \cite{Young-Sym2-LS} may still be adapted  to our setting---albeit  more involved in the dyadic case $\RN (\frp) = 2$  as now $ v_{\frp} (2) = 2 $. 
\end{remark}

By \eqref{4eq: mult}, \eqref{4eq: trivial}--\eqref{4eq: semi-primitive, 2}, one may deduce the following corollary.

\begin{cor}\label{cor: average}
Let $\gamma $ be real. 	We have
	\begin{equation}\label{4eq: average chi0}
		\sum_{\RN (\frc) \leqslant C} \RN (\frc)^{\gamma} \big|  \widehat{F} (\vchi_0)  \big| \Lt  C^{\, \max \{ 1 + \gamma, 0 \}+\vepsilon}, 
	\end{equation}
	 \begin{equation}\label{4eq: average sp}
	 \sum_{\RN (\frc) \leqslant C }	\RN (\frc)^{\gamma} \mathop{\sumsp}_{\vchi (\mathrm{mod}\, \frc)}  	\big|\widehat{F} (\vchi) \big| \Lt C^{\, \max \{ 1 + \gamma, 0 \}+\vepsilon} , 
	 \end{equation}
 where the {\scriptsize \ding{73}}   restricts to semi-primitive $\vchi$. 
\end{cor}

\subsection{Archimedean Mellin Inversion} 

Next, we have the following Mellin representation of $e [z]$ as in \cite[\S 12.2]{Qi-GL(3)} or \cite[\S 2]{Qi-Spectral-LS} (\eqref{2eq: bound xi} is slightly weaker for the convenience of application).  

\begin{lem}\label{lem: Mellin xi}
	Suppose that $ X/\sqrt{2} < |z| < 2 X $.    Then for any $A \geqslant 0$ we have
\begin{align}\label{4eq: Mellin}
	e [z] = \int \!\!\! \int_{\boldsymbol{\EuScript{A}}}  \xi_X  (t      , p          )  {\vlambda_{it      ,  p          } (z)} \mathrm{d}\mu (t      , p          )  , 
\end{align}
 with 
\begin{equation}\label{2eq: bound xi}
	\xi_X (t      , p          ) \Lt \left\{ 
	\begin{split}
		&\displaystyle \frac {\log (2+X)} {1+X}, & & \text{ if } |t      |+|p          | \Lt  1+ X,  \\
		&\displaystyle   {(|t      |+|p          |)^{-A}}, & & \text{  otherwise.  }
	\end{split}
	\right.
\end{equation}

\end{lem}

\subsection{Proof of Proposition \ref{prop: quad form}} 


In order to separate the variables $m$ and $n$, we use the Mellin inversion for  $F (w) $ as  in \eqref{2eq: Mellin inversion}  and for $e [z] $ as in \eqref{4eq: Mellin}. We infer that 
\begin{align*}
F^{\gamma} (d, \theta,   C, M, N) = \sum_{  \RN (c) \sim C}  \!  {\RN (c)^{\gamma}}  \!  \sum_{\vchi (\mathrm{mod}\, c)} \! \widehat{F} (\vchi) \vchi (d) \!  \int \!\!\! \int_{\boldsymbol{\EuScript{A}}}  { 	\xi_{X} (t, p          ) }   {\vlambda}_{it      ,  p          } (\theta/c) B_{it, p          } (\vchi  ) \mathrm{d}\mu  (t, p          ) ,
\end{align*}
where 
$$X = \frac {\sqrt{MN} |\theta| }  {\sqrt{C}}, $$  
and 
\begin{align*}
B_{it, p          } (\vchi  ) =	\bigg( \sum_{\RN (m) \sim M} a_{(m)}  \vchi (m) \vlambda_{it      ,  p          } (m) \bigg) \bigg(  \sum_{\RN (n) \sim N}       \overline{b}_{(n)}  \vchi (n) \vlambda_{it      ,  p          } (n) \bigg). 
\end{align*}
Next we factor $\vchi = \vchi_0 \vchi_1 \vchi_2$ and $\frc = \frc_0 \frc_1 \frc_2$, with $\frc_0$,  $\frc_1$, and $ \frc_2$ relatively prime to each other, characterized by the assumption that $\vchi_0$ is trivial, $\vchi_1$ is primitive, and  $\vchi_2$ is semi-primitive. By \eqref{4eq: mult} we have 
\begin{align*}
	 \big| \widehat{F} (\vchi) \big| =   \big| \widehat{F} (\vchi_0)  \big| \big| \widehat{F} (\vchi_1) \big| \big| \widehat{F} (\vchi_2) \big|. 
\end{align*} 
Consequently, 
\begin{equation*}
	\begin{split}
		  F^{\gamma} (d, \theta,  &  C, M, N)    \Lt  {C^{\gamma}}    \mathop{\sum \sum}_{ \RN (\frc_0 \frc_2) \Lt C }    \big| \widehat{F} (\vchi_0)  \big|  \ \  \sumsp_{ \vchi_2 (\mathrm{mod}\, \frc_2)   }  \big| \widehat{F} (\vchi_2)\big|   \\
	 &\cdot  \sum_{\RN(\frc_1) \Lt C / \RN (\frc_0 \frc_2) } \, \sumx_{ \vchi_1 (\mathrm{mod}\, \frc_1)   }  \big| \widehat{F} (\vchi_1)\big| \int \!\!\! \int_{\boldsymbol{\EuScript{A}}}      |\xi_{X} (t, p          )  |    |B_{it, p          } (\vchi_0 \vchi_1 \vchi_2  ) | \mathrm{d}\mu  (t, p          ) . 
	\end{split}
\end{equation*}
Our last step will be to bound the second line by Cauchy--Schwarz and the hybrid large sieve in Lemma \ref{lem: hybrid}, with the aid of Lemmas \ref{lem: Mellin F (chi)} (1) and \ref{lem: Mellin xi}, and then,  in a trivial manner,  average as in the first line by Corollary \ref{cor: average}.  
More explicitly, it follows from  Cauchy--Schwarz, \eqref{2eq: hybrid} (absorb $\vchi_0 \vchi_2$ into $\boldsymbol{a}$ or $\overline{\boldsymbol{b}}$), \eqref{4eq: mult}, \eqref{4eq: primitive} (indeed $ | \widehat{F} (\vchi_1)  | \leqslant \RN(\frc_1) / \varphi (\frc_1) $), and \eqref{2eq: bound xi} (choose $A = 3$, say) that the second line is bounded by 
\begin{align*}
	     \bigg ( {\frac {C^2 (1+X)^2 } {\RN (\frc_0 \frc_2)^2} +  {M} } \bigg )^{1/2} \bigg (  {\frac {C^2 (1+X)^2 } {\RN (\frc_0 \frc_2)^2} +  {N} } \bigg )^{1/2}   \frac 1 {1+X} (C (1+X))^{\vepsilon} \|\boldsymbol{a}_{M} \|_2  \|\boldsymbol{b}_{N} \|_2  , 
\end{align*}
or, in view of $K = C (1+X)$, by
\begin{align*}
	\bigg( {\frac {   \sqrt{K} } {\RN (\frc_0 \frc_2)} +  \sqrt{\frac{M}{K}} } \bigg) \bigg(   {\frac {  \sqrt{K} } {\RN (\frc_0 \frc_2) } +  \sqrt{\frac{N}{K}} } \bigg)  C  K^{\vepsilon} \|\boldsymbol{a}_{M} \|_2  \|\boldsymbol{b}_{N} \|_2, 
\end{align*}
and the proof is concluded by applying \eqref{4eq: average chi0} and \eqref{4eq: average sp} to the resulting sums over $\frc_0$ and $\frc_2$ respectively. 


\section{Spectral Kuznetsov Trace Formula for \texorpdfstring{$\mathrm{PSL}_2 (\mathbb {Z} [i])$}{PSL\unichar{"2082}(Z[i])}}

Our main tool  is the spectral Kuznetsov trace formula of Bruggeman and Motohashi over the Gaussian field $\mathrm{F}$ \cite[Theorem 10.1]{B-Mo}. Similar to \cite[Lemma 5]{Qi-Spectral-LS}, it is formulated here in our set of notation.

For $\vnu \in \mathbb {C}$, $  p           \in \mathbb {Z}  $, and $z \in \mathbb {C} \smallsetminus \{0\}$,  define 
\begin{equation}\label{0def: J mu m (z)}
	J_{\vnu ,   p          } (z) = J_{\vnu + p           }    (z)   J_{\vnu -  p            }    ({\bar z} ) ,
\end{equation} 
\begin{equation}\label{0eq: defn of Bessel}
	\boldsymbol{J}_{ \vnu,   p          } (z)  =   \frac {2\pi^2} {\sin (\pi \vnu)}  ( J_{-\vnu,\shskip  - p          } (   z) - J_{\vnu,   p          } (    z)  ) ,   
\end{equation}
where $ J_{\nu} (z) $ is the Bessel function of the first kind  (see \cite{Watson}). Note that $ \boldsymbol{J}_{ \vnu,   p          } (z) $ is an even function. 

\begin{lem}\label{lem: Kuznetsov}
	Let $h (t      , p          )$ be an even function on $\boldsymbol{\EuScript{A}}$   that admits an entire analytic continuation $ h (t       + i \sigma, p          ) $ so that it  decays rapidly in both $ t      $ and $ p           $,  uniformly for $\sigma$ on   bounded intervals.  
	For $m , n \in \CaloO \smallsetminus \{0\}$, we have the identity{\rm:} 
	\begin{equation}\label{1eq: Kuznetsov} 
		\begin{aligned}
			\sum_{\uppii \shskip \in \Pi_{c} }       \frac {  \lambda_{\uppii} ( m )     \lambda_{\uppii}  ( n )} {L(1, \mathrm{Sym}^2 \uppii )}   & h  ( t      _{\uppii}, p          _{\uppii} ) +   \frac 1 {\pi}   \int \!\!\! \int_{\boldsymbol{\EuScript{A}}}  
			\frac {\tau_{i t      , p          } (m )  \tau_{i t      , p          } ( n )} {|\zeta (1+2i t      , 2 p          )|^2} h ( 2t      , 4 p           )  \shskip   \mathrm{d} \mu (t      , p          )  \\
			& =        \frac {1} { 8\pi^3 } \, \SDH \sum_{\epsilon \shskip \in \CaloO^{\times  2} \! } \delta_{m,  \epsilon n}  + \frac 1 {32 \pi^3 }    \sum_{c  \shskip \in   \CaloO \smallsetminus \{0\} } \frac {S  (  m ,    n  ; c  ) } { \RN(c) } \SDH  \bigg( \frac { 2\pi \sqrt{  m n} } {    c     }   \bigg),
		\end{aligned}
	\end{equation}
	where  $\delta_{m,  n}$ is the Kronecker $\delta$ symbol, $ \SDH  $ and $ \SDH (z)$ are the Plancherel and Bessel integrals defined by  
	\begin{align}\label{1eq: defn Bessel integral}
		\SDH  = \!  \int \!\!\! \int_{\boldsymbol{\EuScript{A}}}   h (t      , p          )  (t      ^2 + p          ^2 )  \mathrm{d} \mu (t      , p          ), \quad \SDH (z) = \!  \int \!\!\! \int_{\boldsymbol{\EuScript{A}}}   h (t      , p          )  \boldsymbol{J}_{  i t      , p          } ( z )  (t      ^2 + p          ^2 )  \mathrm{d} \mu (t      , p          )  . 
	\end{align} 
\end{lem}

\section{Formulae of the Bessel Integral}

Let    
\begin{equation}\label{4eq: choice of h} 
	h( t,  p           ) = \exp \left ( - \Big(  \frac {t } {T} \Big)^2 - \Big(  \frac {p           } {P } \Big)^2 \right ) .
\end{equation}
As recorded below in Lemma \ref{lem: H(z)}, two formulae for $$\SDH (z) =    \int \!\!\! \int_{\boldsymbol{\EuScript{A}}}   h (t      , p          )  \boldsymbol{J}_{  i t      , p          } ( z )  (t      ^2 + p          ^2 )  \mathrm{d} \mu (t      , p          ) $$ were proved in \cite[\S 4]{Qi-Spectral-LS}---we may transform from one to another by partial integration. 

\begin{lem}\label{lem: H(z)}
Let $\boldsymbol{\hat{\EuScript{A}}}_+ = (-\infty, \infty)  \times [0, \pi / 2) $.	We have 
	\begin{align}\label{4eq: integral H(z)}
		\SDH (z) =   - 2 \int \!\!\! \int_{ {\boldsymbol{\hat{\EuScript{A}}}_+}  } & \cos ( 2 \mathrm{Re} (z \mathrm{trh}(r, \omega)) )    \big(k''(r) \theta (\omega) + k (r) \theta'' (\omega)\big) \mathrm{d} r \shskip \mathrm{d} \omega   , \\
		\label{4eq: integral H(z), 2}
		\SDH (z) =    8  |  z|^2   \int \!\!\! \int_{ \boldsymbol{\hat{\EuScript{A}}}_+  } & \cos (2 \mathrm{Re} (z \mathrm{trh}(r, \omega)) )  (\sinh^2 r  +  \sin^2 \omega  ) k  (  r) \theta (  \omega)  \mathrm{d} r \shskip \mathrm{d} \omega ,
	\end{align} 
	in which
	\begin{align}\label{10eq: trh function, 0}
		\mathrm{trh} (r, \omega) =    \cosh r \cos \omega + i \sinh r \sin \omega,
	\end{align}  
	\begin{align}\label{4eq: defn of h}
		k  (r ) =  \sqrt{\pi} T \exp \big( \!  -  (T r )^2   \big), \qquad 	
		\theta (\omega) =   \sqrt{\pi} P \sum_{q = -\infty}^{\infty}  \exp \big( \!   - (P (\omega + \pi q) )^2 \big)    . 
	\end{align} 
\end{lem}

\begin{remark}
	For our convenience, we have replaced the integral domain $    \boldsymbol{\hat{\EuScript{A}}} = (-\infty, \infty)  \times [0, \pi)$ therein by  $\boldsymbol{\hat{\EuScript{A}}}_+ = (-\infty, \infty)  \times [0, \pi / 2) $, since the integrand is invariant under the changes $r \rightarrow -r$ and $\omega \rightarrow \pi - \omega$. This allows us to avoid the subtle issue that $ \mathrm{trh}(r, 0) = \cosh r  $ and $\mathrm{trh}(r, \pi) = - \cosh r$ have different signs. 
\end{remark}


\section{Proof of Theorem \ref{thm: sym2 large sieve}}


As remarked at the end of \S \ref{sec: main results}, we shall assume that   $ {a}_{\frn} = 0$ outside the dyadic range $ N < \RN (\frn) \leqslant 2N$. Let us assume that $T$ and $P$ are at least $(T PN)^{\vepsilon}$---this may always be ensured by enlarging $T$ or $P$ by $(T PN)^{\vepsilon}$. 

The primary object of study in this section is the sum $\SC + \SE$, defined by   
\begin{align}   \label{7eq: smoothed sum, C}
\SC   =	& \sum_{\uppii \shskip \in \Pi_{c}  }  \omega_{\uppii} { h  ( t      _{\uppii},  p          _{\uppii} ) } \bigg| \sum_{ n  }  a_{(n)} \lambda_{\uppii} (n^2)   \bigg|^2 , \\
\label{7eq: smoothed sum, E}
\SE   =	 \frac 1 {\pi}   	\int \!\!\! \int_{ \boldsymbol{\EuScript{A}} } & \omega (t, p          )   {h (2t      , 4 p          )}   \bigg| \sum_{ n }  a_{(n)} \tau_{i t      , p          } (n^2)   \bigg|^2 \! \mathrm{d}\mu (t      , p           ) . 
\end{align}
In view of our choice of $h (t, p          )$ as in \eqref{4eq: choice of h}, it is clear that $ \SC^{(2)} (T, P, N) \allowbreak + \SE^{(2)} (T, P, N) $  (as defined by \eqref{1eq: defn of C} and \eqref{2eq: defn of E}) is  dominated by $\SC + \SE$.  

\subsection{Application of Kuznetsov}

We start by opening the squares in \eqref{7eq: smoothed sum, C} and \eqref{7eq: smoothed sum, E}, and applying the Kuznetsov trace formula  \eqref{1eq: Kuznetsov} 
so that  \begin{align}
	 \SC + \SE = \SD + \SF,  
\end{align}
where $ \SD$ is the diagonal contribution
\begin{align}
	\SD = O \bigg( T P \big(T^2 + P^2\big)      \sum_{ n } |a_{(n)}|^2 \bigg), 
\end{align}
while $\SF$ is the off-diagonal Kloosterman--Bessel sum 
\begin{align}
	\SF = \frac 1 {32 \pi^3 }     \mathop{\sum\sum}_{m, n}   a_{(m)} \overline{a}_{(n)} \sum_{ c } \frac {S  (  m^2 ,    n^2  ;   c  ) } { \RN(c) } \SDH  \bigg( \frac { 2\pi  {    m n} } {    c     }   \bigg) . 
\end{align} 

Next we  express the Bessel integral $ \SDH (2\pi   mn/  c) $ by the formulae in \eqref{4eq: integral H(z)} and \eqref{4eq: integral H(z), 2} and extract the exponential factor $e [2   mn/c]$ to join the Kloosterman sum $ S  (   m^2 ,    n^2  ; c  ) $ via the identity
 \begin{align*}
 	\cos \Big( 4\pi \mathrm{Re} \Big(\frac { mn} {c} \mathrm{trh}(r, \omega) \Big) \Big) = \frac 1 2 \sum_{\pm}  e \Big[ \pm \frac {2  mn} {c} \Big] e \Big[ \pm \frac {2  mn} {c}  ( \mathrm{trh}(r, \omega)-1)\Big]   .
 \end{align*}
Note that the $\pm$-sum  may be absorbed into the $c$-sum. It follows that 
\begin{align}\label{7eq: F, 1}
	\SF   = - \frac 1 {16 \pi^3 }     \int \!\!\! \int_{ {\boldsymbol{\hat{\EuScript{A}}}_+}  }   \SF_0  (r, \omega) \big(k''(r) \theta (\omega) + k (r) \theta'' (\omega)\big) \mathrm{d} r \shskip \mathrm{d} \omega ,  
\end{align}
\begin{align}\label{7eq: F, 2} 
	\SF  = \frac 1 {  \pi}   \int \!\!\! \int_{ {\boldsymbol{\hat{\EuScript{A}}}_+}  }   \SF_1  (r, \omega)  (\sinh^2 r  +  \sin^2 \omega  ) k  (  r) \theta (  \omega)  \mathrm{d} r \shskip \mathrm{d} \omega  ,   
\end{align}
where, for $\delta =0$ or $1$, 
\begin{align}\label{6eq: Fd}
	\SF_{\delta} (r, \omega) = \mathop{\sum\sum}_{m, n}  \RN(mn)^{\delta}  a_{(m)} \overline{a}_{(n)} \cdot D_{\delta} (m, n; r, \omega) , 
\end{align} 
with
\begin{align}\label{6eq: Dd}
	D_{\delta} (m, n; r, \omega) =   \sum_{ c } \frac {S  (   m^2 ,    n^2  ; c  ) e[2  mn/c] } { \RN(c)^{\delta + 1} } e \Big[ \frac {  mn} {c}  \psi (r, \omega) \Big] , 
\end{align} 
\begin{align}
	 \psi (r, \omega) = 2(\mathrm{trh}(r, \omega) -1). 
\end{align}

\subsection{Reduction} 

Following  the proof of Lemma 3.1 in  \cite{Young-Sym2-LS}, we would like to reformulate $D_{\delta} (m, n; r, \omega) $.  To this end,   let us invoke the Selberg identity
 \begin{align*}
 	S (m^2, n^2; c) = \sum_{\mathfrak{d} | (m^2, n^2, c)}  \RN (d) S \Big( \Big ( \frac{mn} {d } \Big )^2, 1; \frac c d \Big)  ,
 \end{align*}
with $\mathfrak{d} = (d)$,  and the formula 
\begin{equation*}
	S (w^2, 1; c g) = \left\{ \begin{split}
	&	0, & &  \text{ if } (c, g) \neq (1), \\
	& \mu (g) S \big((w/g)^2, 1; c \big), & &  \text{ if } (c, g) = (1),
	\end{split}\right.
\end{equation*}
for $ w \equiv 0\, (\mathrm{mod} \, g)$. 
Thus,  if we set $g = (mn/d, c/d)$ and change   $c$ (as in \eqref{6eq: Dd}) into $ d g \cdot c $,  some direct calculations  yield 
\begin{align}\label{6eq: D}
	D_{\delta} (m, n; r, \omega) =    \mathop{\mathop{\sum\sum}_{\mathfrak{d}| (m^2, n^2)}   }_{\mathfrak{d} \mathfrak{g} | (mn)}   \frac{\mu (\mathfrak{g}) } {\RN(\mathfrak{d})^{\delta } \RN(\mathfrak{g})^{\delta+1}}  F_{\delta} \bigg( \frac {  mn} { d g}  ; \psi (r, \omega) \bigg) , 
\end{align}
for $\delta = 0$ or $1$, where by our convention $\mathfrak{d} = (d)$, $\mathfrak{g} = (g)$, and we define 
\begin{align}\label{7eq: defn of F}
F_{\delta} (w ; \psi) =  \sum_{ (c, w) = (1)}   \frac { F   (  w ; c ) e   [ w  \psi  / c  ]  } {\RN (c  )^{\delta+1}} , \qquad F  (w; c) = S (w^2, 1; c) e  [2 w / c ]. 
\end{align}
Note that if we decompose $\mathfrak{d} = \mathfrak{d}_1 \mathfrak{d}_2^2$ with $\mathfrak{d}_1$ square-free, then the condition $ \mathfrak{d} | (m^2, n^2) $ is equivalent to $ \mathfrak{d}_1 \mathfrak{d}_2 | (m, n) $. Consequently, after inserting \eqref{6eq: D} into \eqref{6eq: Fd},  we may use  the substitutions $ m \rightarrow    d_1 d_2 m$ and $n \rightarrow    d_1 d_2 n $   to transform $ \SF_{\delta} (r, \omega) $ into 
\begin{equation*}
	  \begin{split}
	 \SF_{\delta} (r, \omega)=  \! \!  \mathop{\sumf \! \sum}_{ \mathfrak{d}_1,  \mathfrak{d}_2  }   {\RN ( \mathfrak{d}_1  )^{\delta } } \mathop{\sum\sum}_{m, n} \RN (mn)^{\delta } a_{m\mathfrak{d}_1 \mathfrak{d}_2} \overline{a}_{n\mathfrak{d}_1 \mathfrak{d}_2}    \sumf_{\mathfrak{g} | mn\mathfrak{d}_1} \! \frac {\mu (\mathfrak{g})} {\RN ( \mathfrak{g} )^{\delta+1}} F_{\delta} \bigg( \! \frac { d_1 mn } {   g}  ; \psi (r, \omega) \! \bigg) ,
	  \end{split}
\end{equation*} 
where the superscript $\flat$ indicates summation over square-free ideals.  Further, we decompose $\mathfrak{d}_1 = \mathfrak{d}_{\flat} \mathfrak{g}_{1} $ and $\mathfrak{g} = \mathfrak{g}_1 \mathfrak{h} \mathfrak{k}$ by  
\begin{align*}
	 \mathfrak{g}_{1}  = (\mathfrak{g}, \mathfrak{d}_1 ), \qquad \mathfrak{h} = (\mathfrak{g}\mathfrak{g}_{1} ^{-1}, m),  
\end{align*} 
and make the changes $ m \rightarrow  h m $ and $n \rightarrow   k n$. 
 It follows that 
 \begin{equation}\label{7eq: final sum} 
 	\begin{split}
 	 \SF_{\delta} (r, \omega) =     \mathop{\sumf \sumf }_{(\mathfrak{d}_{\flat}, \mathfrak{g}_1) = (1)} \sum_{ \mathfrak{d}_2}   {\RN ( \mathfrak{d}_{\flat}  )^{\delta} }    \mathop{\mathop{\sumf  \sumf}_{  (\mathfrak{h} ,  \mathfrak{k}) = (1) } }_{(\mathfrak{h k}, \mathfrak{d}_{\flat}) = (1) }  \frac { \mu (    \mathfrak{g}_1 \mathfrak{h} \mathfrak{k} )} {\RN (\mathfrak{g}_1      \mathfrak{h} \mathfrak{k} ) }   \mathop{\mathop{\sum\sum}_{m, n}}_{(m, \mathfrak{k}) = (1)}  a^{\delta, \mathfrak{h}}_{(m)} \bar{a}{}^{\delta, \mathfrak{k}}_{(n)}   F_{\delta}   \big( d_{\flat} mn  ; \psi (r, \omega)  \big) , 
 	\end{split}
 \end{equation} 
with the definitions
\begin{align}\label{7eq: new sequences}
a^{\delta, \mathfrak{h}}_{\frm} = \RN (\frm)^{\delta }	a_{\frm \mathfrak{h} \mathfrak{d}_{\flat} \mathfrak{g}_1 \mathfrak{d}_2}, \qquad a^{\delta, \mathfrak{k}}_{ \frn } = \RN (\frn)^{\delta } {a}_{\frn \mathfrak{k} \mathfrak{d}_{\flat} \mathfrak{g}_1  \mathfrak{d}_2}  . 
\end{align}


\subsection{Truncation and Partition}\label{sec: truncation}
Next, some cleanup by truncation and partition is needed  prior to the application of Proposition \ref{prop: quad form}  to the inner sum in  the expression of $ \SF_{\delta} (r, \omega)  $ as in  \eqref{7eq: final sum}. 

By the Weil bound (see \cite[Theorem 10]{BM-Kloostermann})
\begin{align*}
	S (m, n; c)   \Lt \tau (c) \sqrt{\RN(m, n, c)} \sqrt{\RN(c)}, 
\end{align*} 
we infer that if  the $c$-sum is truncated at $ \RN (c) = N^2 $, then the tail of the sum $ \SF_1 (r, \omega) $ as in \eqref{6eq: Fd} may be trivially bounded by $ O \big( N^{2+\vepsilon} 	\| \boldsymbol{a}_N \|_2^2 \big) $, which is satisfactory as the integration in \eqref{7eq: F, 2} produces an extra factor $O \big(1/T^2+1/P^2 \big)$.  

For dyadic $C < N^2$, let $F_{\delta} (w ; \psi; C)$ 
denote the partial sum of $ F_{\delta} (w ; \psi) $ (see \eqref{7eq: defn of F}) over $ C < \RN (c) \leqslant 2C  $.  Let $ \SF_{\delta} (r, \omega; C)   $ and $\SF (C)$ be defined accordingly (see \eqref{7eq: F, 1}, \eqref{7eq: F, 2}, and \eqref{7eq: final sum}). In particular, write the inner sum of $ \SF_{\delta} (r, \omega; C)   $ (see \eqref{7eq: final sum}) as
\begin{align}\label{7eq: defn Fn(C)}
\SF_{\delta}^{\mathfrak{h}, \mathfrak{k} } (r, \omega; C) =	\mathop{\mathop{\sum\sum}_{m, n}}_{(m, \mathfrak{k}) = (1)}  a^{\delta, \mathfrak{h}}_{(m)} \bar{a}{}^{\delta, \mathfrak{k}}_{(n)}   F_{\delta}   \big( d_{\flat} mn  ; \psi (r, \omega) ; C \big). 
\end{align}

By the exponential decay of $ k (r) $ and $\theta (\omega)$ as in \eqref{4eq: defn of h}, we may truncate the integral in \eqref{7eq: F, 1} or \eqref{7eq: F, 2} at $|r| = T^{\vepsilon}/ T$ and $   \omega = P^{\vepsilon}/ P$, at the cost of a negligibly small error. Define  
\begin{align*}
	\boldsymbol{\hat{\EuScript{A}}}^{\natural}_+  = \big\{ (r, \omega) :|r| \leqslant T^{\vepsilon}/ T,   \  \omega \leqslant P^{\vepsilon}/ P  \big\}. 
\end{align*}

Further, in order to facilitate our later analysis, we partition the $(r, \omega)$-integral according to the value of 
\begin{align*}
	|\psi (r, \omega)| = 2 (\cosh r - \cos \omega) = 4  ( \sinh^2 (r/2) + \sin^2 (\omega/2)  ). 
\end{align*} 
To this end, for small $\phi \Lt 1$ we introduce
\begin{align*}
	\boldsymbol{\hat{\EuScript{A}}}_{\mathrm{o}} (\phi) = \big\{ (r, \omega) : r^2 + \omega^2 \leqslant \phi \big\}, \quad   \boldsymbol{\hat{\EuScript{A}}}_+  (\phi) = \big\{ (r, \omega) : \phi < r^2 + \omega^2 \leqslant 2 \phi \big\}. 
\end{align*}
By symmetry, let us assume $P \leqslant T$ for simplicity.  Set 
\begin{align*}
	\phi_{\mathrm{o}} = \frac {\sqrt{C}} {N}, \qquad   \phi_{\natural} =   \frac { P^{\vepsilon}}  { P^2 }. 
\end{align*}
Now we partition further the truncated integral domain $ \boldsymbol{\hat{\EuScript{A}}}^{\natural}_+ $   into the union of $   \boldsymbol{\hat{\EuScript{A}}}^{\natural}_{\mathrm{o}}   = \boldsymbol{\hat{\EuScript{A}}}^{\natural}_+ \cap \boldsymbol{\hat{\EuScript{A}}}_{\mathrm{o}} (\phi_{\mathrm{o}})$ and  those $ \boldsymbol{\hat{\EuScript{A}}}^{\natural}_+ (\phi) = \boldsymbol{\hat{\EuScript{A}}}^{\natural}_+ \cap \boldsymbol{\hat{\EuScript{A}}}_+ (\phi) $ with dyadic $ \phi_{\mathrm{o}} \leqslant \phi  \Lt \phi_{\natural} $. 
Accordingly, we let $ \SF_{\mathrm{o}} (C)$ and $\SF (\phi; C)$ denote the contributions from these domains to the integral $\SF (C)$ (see  \eqref{7eq: F, 1} or \eqref{7eq: F, 2}). 

\subsection{Application of Proposition \ref{prop: quad form}} 

Now we are ready to apply Proposition \ref{prop: quad form} to the sum $\SF_{\delta}^{\mathfrak{h}, \mathfrak{k} } (r, \omega; C)$ as defined by \eqref{7eq: new sequences} and \eqref{7eq: defn Fn(C)} (see also \eqref{7eq: defn of F}). Set
\begin{align*}
	N_{\mathfrak{h}} = \frac {N} { \RN (\mathfrak{h} \mathfrak{d}_{\flat} \mathfrak{g}_1 \mathfrak{d}_2  ) }, \qquad N_{\mathfrak{k}} = \frac {N} { \RN (\mathfrak{k} \mathfrak{d}_{\flat} \mathfrak{g}_1 \mathfrak{d}_2  ) }. 
\end{align*} Note that $|\theta| = \sqrt{\RN (\mathfrak{d}_{\flat})}  |\psi (r, \omega)|$ and hence 
\begin{align*}
	K = C +   {\textstyle \sqrt{C N_{\mathfrak{h}} N_{\mathfrak{k}} \RN (\mathfrak{d}_{\flat})}   } |\psi (r, \omega)| = C +  \frac {\sqrt{C} N |\psi (r, \omega)|} {\sqrt{\RN (\mathfrak{hk} \mathfrak{d}_{\flat})} \RN (\mathfrak{g}_1 \mathfrak{d}_2 )}. 
\end{align*} 
So if we were to apply \eqref{1eq: Kloosterman} directly, the expressions would be quite involved with $ \RN (\mathfrak{h}) $, $\RN (\mathfrak{k})$, ..... in the denominators. However, it is quite safe to drop them all, as we shall only need $1 / \RN (\mathfrak{hk})$ while this has already occurred in \eqref{7eq: final sum}. 
By applying  \eqref{1eq: Kloosterman} in this way, for $(r, \omega)$ in $ \boldsymbol{\hat{\EuScript{A}}}_{\mathrm{o}} (\phi_{\mathrm{o}}) $ ($\phi_{\mathrm{o}} = \sqrt{C}/ N$),  we have
\begin{align}\label{7eq: bound for F, 1.1}
	 \SF_{\delta}^{\mathfrak{h}, \mathfrak{k} } (r, \omega; C) \Lt \frac {C+N} {C^{\delta }} 
	 N^{\vepsilon} \big\|\boldsymbol{a}^{\delta, \mathfrak{h}}_{N_{\mathfrak{h}}}\big\|_2 \big\|\boldsymbol{a}^{\delta, \mathfrak{k}}_{N_{\mathfrak{k}}} \big\|_2, 
\end{align}
while for $ (r, \omega) $ in $ \boldsymbol{\hat{\EuScript{A}}}_+  (\phi)$ (recall that $  \phi  \geqslant \sqrt{C}/ N$), we have 
\begin{equation}\label{7eq: bound for F, 1.2} 
		\SF_{\delta}^{\mathfrak{h}, \mathfrak{k} } (r, \omega; C) \Lt \frac 1 {C^{\delta }} \bigg(     {\sqrt{C} N \phi }     + \frac {\sqrt{C}} {\phi}  \bigg) N^{\vepsilon} \big\|\boldsymbol{a}^{\delta, \mathfrak{h}}_{N_{\mathfrak{h}}}\big\|_2 \big\|\boldsymbol{a}^{\delta, \mathfrak{k}}_{N_{\mathfrak{k}}} \big\|_2;  
\end{equation}
here we have used the simple inequalities $ C \sqrt{N_{\mathfrak{h}} N_{\mathfrak{k}}}/ K < N$ and $ C \sqrt{N_{\mathfrak{h}} N_{\mathfrak{k}}}/ K \Lt \sqrt{C  }/ \phi$ (see  \eqref{1eq: Kloosterman}) in these two cases.    

Similar to \eqref{7eq: bound for F, 1.1} and \eqref{7eq: bound for F, 1.2}, we claim that 
\begin{align}\label{7eq: bound for F, 2.1}
	 \SF_{\delta}  (r, \omega; C) \Lt (N^2/C)^{\delta } (C+N) 
	 N^{ \vepsilon}  \|\boldsymbol{a}_N \|_2^2, 
\end{align}
for  $(r, \omega)$ in $ \boldsymbol{\hat{\EuScript{A}}}_{\mathrm{o}} (\phi_{\mathrm{o}}) $, and
\begin{align}\label{7eq: bound for F, 2.2}
	\SF_{\delta}  (r, \omega; C) \Lt {\sqrt{C}} (N^2/C)^{\delta }  (     N \phi      +   1 / {\phi}   ) 
	N^{ \vepsilon}  \|\boldsymbol{a}_N \|_2^2, 
\end{align}
for  $ (r, \omega) $ in $ \boldsymbol{\hat{\EuScript{A}}}_+  (\phi)$. Indeed, in view of \eqref{7eq: final sum}--\eqref{7eq: defn Fn(C)}, \eqref{7eq: bound for F, 2.1} and \eqref{7eq: bound for F, 2.2} may be deduced from \eqref{7eq: bound for F, 1.1} and \eqref{7eq: bound for F, 1.2} by the following
\begin{align*}
	 \mathop{\sum \sum \sum}_{ \mathfrak{d}_{\flat},  \mathfrak{g}_1 , \mathfrak{d}_2 }    {\RN ( \mathfrak{d}_{\flat}  )^{\delta } }    \mathop{\sum  \sum }_{  \mathfrak{h} ,  \mathfrak{k} }   \frac { 1 } {\RN (\mathfrak{g}_1      \mathfrak{h} \mathfrak{k} ) } \bigg( \frac {N^2}  { {\RN (\mathfrak{h k})} \RN (  \mathfrak{d}_{\flat} \mathfrak{g}_1 \mathfrak{d}_2  )^2  } \bigg)^{\delta}  \big\|\boldsymbol{a}^{ 0, \mathfrak{h}}_{N_{\mathfrak{h}}}\big\|_2 \big\|\boldsymbol{a}^{ 0, \mathfrak{k}}_{N_{\mathfrak{k}}} \big\|_2 \Lt N^{ 2\delta + \vepsilon} \|\boldsymbol{a}_N \|_2^2,  
\end{align*}
where we have dropped the co-primality and square-free conditions in \eqref{7eq: final sum}.  
Next, we drop all except $\RN (\mathfrak{hk})$ in the denominator, group $ \mathfrak{d}_{\flat}$,  $\mathfrak{g}_1 $, and $ \mathfrak{d}_2 $ into a new variable $ \mathfrak{t} = \mathfrak{d}_{\flat} \mathfrak{g}_1 \mathfrak{d}_2 $, and use the symmetry in the $(\mathfrak{h}, \mathfrak{k})$-sum to form a square. It follows that the left-hand side  is bounded by 
\begin{align*}
	N^{2\delta} \sum_{\mathfrak{t}} \bigg( \sum_{\mathfrak{k}} \frac 1 {\RN (\mathfrak{k})} \big\|\boldsymbol{a}^{ 0,  \mathfrak{k}}_{N_{\mathfrak{k}}} \big\|_2 \bigg)^2, \qquad \text{($a^{0, \mathfrak{k}}_{ \frn } =   {a}_{\frn \mathfrak{k} \mathfrak{t} }$)}, 
\end{align*} 
and by Cauchy--Schwarz this is clearly bounded by $ O \big( N^{ 2\delta + \vepsilon} \|\boldsymbol{a}_N \|_2^2 \big) $ as desired.

\subsection{Final Estimates} 

In view of \eqref{7eq: F, 1}, \eqref{7eq: F, 2}, and the partitions made in  \S \ref{sec: truncation}, our problem is reduced to estimating the integrals 
\begin{align}\label{7eq: F(phi;C), 1}
	\SF  (\phi; C)  = & - \frac 1 {16 \pi^3 }     \int \!\!\! \int_{ {\boldsymbol{\hat{\EuScript{A}}}^{\natural}_+ (\phi)}  }   \SF_0  (r, \omega; C) \big(k''(r) \theta (\omega) + k (r) \theta'' (\omega)\big) \mathrm{d} r \shskip \mathrm{d} \omega ,  \\
	\label{7eq: F(phi;C), 2} 
	\SF  (\phi; C) & = \frac 1 {  \pi}   \int \!\!\! \int_{ {\boldsymbol{\hat{\EuScript{A}}}^{\natural}_+ (\phi)}  }   \SF_1  (r, \omega; C)  (\sinh^2 r  +  \sin^2 \omega  ) k  (  r) \theta (  \omega)  \mathrm{d} r \shskip \mathrm{d} \omega  ,   
\end{align}
and 
\begin{align}\label{7eq: Fo(C), 1}
	\SF _{\mathrm{o}}   (C)   = & -   \frac 1 {16 \pi^3 }       \int \!\!\! \int_{ {\boldsymbol{\hat{\EuScript{A}}}^{\natural}_{\mathrm{o}}  }  }   \SF_0  (r, \omega; C) \big(k''(r) \theta (\omega) + k (r) \theta'' (\omega)\big) \mathrm{d} r \shskip \mathrm{d} \omega ,   \\
	\label{7eq: Fo(C), 2} 
 	\SF _{\mathrm{o}}  (  C  )   & = \frac 1 {  \pi}    \int \!\!\! \int_{ {\boldsymbol{\hat{\EuScript{A}}}^{\natural}_{\mathrm{o}} }  }   \SF_1  (r, \omega; C)  (\sinh^2 r  +  \sin^2 \omega  ) k  (  r) \theta (  \omega)  \mathrm{d} r \shskip \mathrm{d} \omega .     
\end{align} 

Recall that we assumed $P \leqslant T$, so our aim as in \eqref{1eq: sym2 large sieve, 1} or \eqref{1eq: sym2 P<T} is to prove that  
$ \SF (\phi; C) $ and $  \SF _{\mathrm{o}} (C) $  are both bounded by 
\begin{align*}
	\bigg(  T P N +   \frac {T N^2} {P^3} \bigg) (TN)^{\vepsilon} \|\boldsymbol{a}_N\|_2^2,
\end{align*}
as long as $\sqrt{C} < N$ and $\sqrt{C}/N \leqslant \phi \Lt P^{\vepsilon}/ P^2$.

Note that $$ \textit{Area} \, (\boldsymbol{\hat{\EuScript{A}}}^{\natural}_+ (\phi) ) \Lt \sqrt{\phi} T^{\vepsilon} / T	 . $$ 
For     $ N\sqrt{\phi} /T \leqslant \sqrt{C} < N$, by  \eqref{7eq: bound for F, 2.2} and   \eqref{7eq: F(phi;C), 2}   we get
\begin{align*}
	\SF  (\phi; C) & \Lt \frac {\sqrt{\phi} } {T} T P \phi \frac {N^2}{\sqrt{C}} \bigg( N \phi + \frac {1} {\phi}\bigg) (TN)^{\vepsilon}   \|\boldsymbol{a}_N \|_2^2 	\Lt T P \big( N  + N^2\phi^2 \big) (TN)^{\vepsilon}   \|\boldsymbol{a}_N \|_2^2. 
\end{align*} 
For     $  \sqrt{C} < N\sqrt{\phi} /T $, by   \eqref{7eq: bound for F, 2.2} and  \eqref{7eq: F(phi;C), 1}    we get
\begin{align*}
	\SF  (\phi; C) \Lt \frac {\sqrt{\phi} } {T} T^3 P \sqrt{C} \bigg( N \phi + \frac {1} {\phi}\bigg) (TN)^{\vepsilon}   \|\boldsymbol{a}_N \|_2^2 \Lt   T P \big(  N + N^2\phi^2 \big) (TN)^{\vepsilon}   \|\boldsymbol{a}_N \|_2^2 . 
\end{align*} 
These meet  the desired bound as $\phi \Lt P^{\vepsilon}/ P^2$. 

Note that  
\begin{equation*}
\textit{Area} \, (\boldsymbol{\hat{\EuScript{A}}}^{\natural}_{\mathrm{o}} )	\Lt \left\{ \begin{split} 
	&  T^{\vepsilon} / T P, & & \text{ if } N/P^2 \leqslant \sqrt{C} < N , \\
	& \sqrt[4]{C} T^{\vepsilon} / \sqrt{N}T, & & \text{ if } N/T^2 \leqslant \sqrt{C} < N/P^2 , \\
	& \sqrt{C} / N, & & \text{ if }   \sqrt{C} < N / T^2,
 	\end{split}  \right. 
\end{equation*}   
and that for $(r, \omega)$ on $ \boldsymbol{\hat{\EuScript{A}}}^{\natural}_{\mathrm{o}}  $ 
\begin{equation*}
	  \sinh^2 r  +  \sin^2 \omega	\Lt \left\{ \begin{split} 
		&  P^{\vepsilon} / P^2, & & \text{ if } N/P^2 \leqslant \sqrt{C} < N , \\
		& \sqrt{C}   / N , & & \text{ if } N/T^2 \leqslant \sqrt{C} < N/P^2 .
	\end{split}  \right. 
\end{equation*}  
For     $ N/P^2 \leqslant \sqrt{C} < N$, by       \eqref{7eq: bound for F, 2.1} and  \eqref{7eq: Fo(C), 2} we get
\begin{align*}
	\SF_{\mathrm{o}} (C) \Lt \frac 1 {T P}  \frac {T P} {P^2} \frac {N^2} {C} (C+N) (TN)^{\vepsilon}   \|\boldsymbol{a}_N \|_2^2 \Lt \bigg(  P^2 N + \frac {N^2 } {P^2 }  \bigg) (TN)^{\vepsilon}   \|\boldsymbol{a}_N \|_2^2. 
\end{align*} 
For $N/T^2 \leqslant \sqrt{C} < N/P^2$, by      \eqref{7eq: bound for F, 2.1} and   \eqref{7eq: Fo(C), 2} we get
\begin{align*}
	\SF_{\mathrm{o}} (C) \Lt \frac {\sqrt[4]{C} } {\sqrt{N} T}  \frac {\sqrt{C}T P} {N} \frac {N^2} {C} (C+N) (TN)^{\vepsilon}   \|\boldsymbol{a}_N \|_2^2 \Lt \bigg(T P N + \frac {N^2 } {P^2 }  \bigg) (TN)^{\vepsilon}   \|\boldsymbol{a}_N \|_2^2. 
\end{align*} 
For $ \sqrt{C} < N /T^2 $, by      \eqref{7eq: bound for F, 2.1} and     \eqref{7eq: Fo(C), 1}  we get
\begin{align*}
	\SF_{\mathrm{o}} (C) \Lt \frac {\sqrt{C}} { N } T^3 P (C+N) N^{\vepsilon} \|\boldsymbol{a}_N \|_2^2 \Lt  \bigg( T P N + \frac {P N^2 } {T^3 }  \bigg) N^{\vepsilon}   \|\boldsymbol{a}_N \|_2^2 .
\end{align*}  
It is clear that all of these are satisfactory as $P \leqslant T$.


\section{Proof of Theorem \ref{thm: mean-Lindelof}}\label{sec: mean-Lindelof}

Let $\uppii \in \Pi_{c} $.  For $\mathrm{Re} (s) > 1$, we have 
\begin{align*}
	L (s, \mathrm{Sym}^2 \uppii ) = \zeta_{\mathrm {F}} (2s) \sum_{\frn \subset \CaloO} \frac {\lambda_{\uppii} (\frn^2) } {\RN(\frn)^{s}}.  
\end{align*}
\delete{We may also write
\begin{align*}
	L (s, \mathrm{Sym}^2 \uppii ) = \sum_{\frn \subset \CaloO} \frac {c_{  \uppii}  (\frn) } {\RN(\frn)^{s} }  ,  \qquad c_{  \uppii}  (\frn) = \sum_{ \mathfrak{l}^2 \mathfrak{m} = \frn } \lambda_{\uppii} (\frm^2). 
\end{align*}
}
The gamma factor of $\mathrm{Sym}^2 \uppii $ is equal to $ (2\pi)^{-3 s} \gamma (s, t    _{\uppii}, p              _{\uppii} ) $, with 
\begin{align*}
	\gamma (s, t    , p               ) = \Gamma (s) \Gamma (s+it     +|p              |) \Gamma (s-it     + |p              |). 
\end{align*}
It is known   that $ L (s, \mathrm{Sym}^2 \uppii ) $ is entire (see \cite{Shimura-Sym2,GJ-GL(2)-GL(3)}) and satisfies the functional equation
\begin{align*}
	\Lambda (s, \mathrm{Sym}^2 \uppii ) = \Lambda (1-s, \mathrm{Sym}^2 \uppii ), 
\end{align*}
where $ \Lambda (s, \mathrm{Sym}^2 \uppii ) = \pi^{-3s} \gamma (s, t    _{\uppii}, p              _{\uppii} ) L (s, \mathrm{Sym}^2 \uppii) $. 
Define 
\begin{align*}
	\RC (s, t, p) = |s|  \, | s +it+|p|| \,  |s-it+|p||. 
\end{align*}

Now let $ \uppii \in \Pi_{c} (T, P) $. Subsequently, we always assume that $\mathrm{Re} (s) = 1/2$ and that $T, P \Gt 1$ are large.  By `negligibly small' we mean $ O_A ( (T^2+P^2)^{-A} )$ for any $A > 0$. Moreover, it will be harmless to remove the harmonic weight $\omega_{\uppii}$ in view of  
\begin{align*}
 {L (1, \mathrm{Sym}^2 \uppii)} \Lt   | i t_{\uppii} + p_{\uppii}  |^{ \vepsilon},  
\end{align*}
by \cite[Theorem 1]{Molteni-L(1)} or \cite[Theorem 2]{Li-L(1)} (as pointed out by a referee, the application of the former requires the automorphy of $\mathrm{Sym}^2 \uppii$ due to Gelbart--Jacquet \cite{GJ-GL(2)-GL(3)}).

In the case $ |s| \geqslant \sqrt{T^2+P^2} $, we have 
\begin{align*}
	\RC (s, t_{\uppii}, p_{\uppii}) \Lt |s|^3 , 
\end{align*}
so the trivial convexity bound for $ L (s, \mathrm{Sym}^2 \uppii) $ yields 
\begin{align*}
	 \sum_{\uppii \shskip \in \Pi_{c} (T, P) }   |L(s , \mathrm{Sym}^2 \uppii )|^2 \Lt |s|^{3+\vepsilon} TP \big(T^2+P^2\big) , 
\end{align*}
as desired.

Next, consider the case $|s| < \sqrt{T^2+P^2}$ so that 
\begin{align*}
	 \RC (s, t_{\uppii}, p_{\uppii}) \Lt |s|  \big(T^2+P^2\big). 
\end{align*}
By the approximate functional equation in \cite[Theorem 5.3]{IK},  we express $L (s, \mathrm{Sym}^2 \uppii )$ by the sum
\begin{align*} 
	   \sum_{ \frn \subset \CaloO }  \frac {   \lambda_{\uppii}  (\frn^2 )   } { \RN (\frn)^s  }     V_{s}     ( \pi^3  \RN(\frn) ; t_{\uppii}, p_{\uppii}  ) + \epsilon (s, t_{\uppii}, p_{\uppii}) \sum_{ \frn \subset \CaloO }  \frac {   \lambda_{\uppii}  (\frn^2 )   } { \RN (\frn)^{1-s}  }     V_{1-s}     ( \pi^3  \RN(\frn)  ; t_{\uppii}, p_{\uppii} )  , 
\end{align*}
where $$ \epsilon (s, t, p) = \pi^{3 (2s-1)} \frac {\gamma (1-s, t, p) } {\gamma (s, t, p) } $$ has  norm unity (for  $\mathrm{Re} (s) = 1/2$).  \cite[Proposition 5.4]{IK} implies that  one may effectively restrict the sums above to the range  $ \RN (\frn) \leqslant (|s|(T^2+P^2))^{1+\vepsilon} $ at the cost of a negligibly small error. In order to separate $ \RN (\frn) $ from $( t_{\uppii}, p_{\uppii} )$, we use  the following expression  due to Blomer \cite[Lemma 1]{Blomer}: 
\begin{align*}
	V_{s}  (y; t      ,  p)  = \frac 1 {2   \pi i }   \int_{ \vepsilon - i U}^{\vepsilon + i U}    y^{ - v}  \zeta_{\mathrm {F}} (2s+2v) \frac {\gamma (s+ v, t      , p)} {\gamma (s, t      , p)} \exp({v^2})     \frac {\mathrm{d}v} {v} +  {O_{\vepsilon}} \bigg( \frac { (T^2+P^2)^{  \vepsilon} } {y^{ \vepsilon} \exp ({U^2 /2 }) } \bigg) .
\end{align*}
 The error term above is again negligibly small if we choose $ U = \log ( T^2+P^2)$. By the Stirling formula, for any $v$ on the integral contour and $ (t      , p) $ in the rectangle $  \boldsymbol{\EuScript{A}} (T, P) $,
\begin{align*}
	\frac {\gamma (s+ v, t      , p)} {\gamma (s, t      , p)}   = O_{\vepsilon} \big( (T^2+P^2)^{\vepsilon}\big). 
\end{align*} 
At any rate, by Cauchy--Schwarz and an application of Theorem  \ref{thm: sym2 large sieve} inside the integral, one may easily derive the bound 
\begin{align*}
	 \sum_{\uppii \shskip \in \Pi_{c} (T, P) }   |L(s , \mathrm{Sym}^2 \uppii )|^2    \Lt  |s| TP \big(T^2 + P^2\big)^{1+\vepsilon} + |s|^2 \frac  { (T^2+P^2 )^{4+\vepsilon}} {T^3 P^3}. 
\end{align*}
By symmetry, let us consider the case $P \leqslant T$ so that
\begin{align*}
	\sum_{\uppii \shskip \in \Pi_{c} (T, P) }   |L(s , \mathrm{Sym}^2 \uppii )|^2    \Lt  |s|  P T^{3+\vepsilon} + |s|^2 \frac  { T^{5+\vepsilon}} { P^3}. 
\end{align*}
Since the left-hand side is non-decreasing in $P$, we may replace $P$ in the right-hand side by $P+\sqrt{T}$, and hence
 \begin{align*}
 	\sum_{\uppii \shskip \in \Pi_{c} (T, P) }   |L(s , \mathrm{Sym}^2 \uppii )|^{2}    \Lt  |s|^{2}  (P+\sqrt{T}) T^{3+\vepsilon},
 \end{align*}
 as desired. 
 
\newcommand{\etalchar}[1]{$^{#1}$}
\def\cprime{$'$}

\end{document}